\documentclass[11pt, a4paper, english]{article}
\usepackage{graphicx} 

\title{Statistical guarantees for signature kernel ridge regression}
\author{Christian Bayer, Jakob Kellermann}
\date{\today}

\usepackage{graphicx} 
\usepackage{xcolor}
\usepackage{amsmath}
\usepackage{algpseudocode}
\usepackage{algorithm}
\usepackage{bm}
\usepackage{shuffle}
\usepackage{booktabs} 
\usepackage{adjustbox}

\usepackage[T1]{fontenc}
\usepackage[utf8]{inputenc}
\usepackage[hidelinks]{hyperref}
\usepackage[style=alphabetic,maxnames=10]{biblatex}
\usepackage{mathtools}

\def\mb{\ensuremath\mathbb}
\def\mc{\ensuremath\mathcal}
\usepackage{listings}
\usepackage{amsmath}
\usepackage{amssymb}
\usepackage{amsthm}
\usepackage{bbm}
\usepackage{mathrsfs}

\usepackage{todonotes}

\newcommand{\bx}{\mathbf{x}}
\newcommand{\bX}{\mathbf{X}}
\newcommand{\mch}{\mathcal{H}}
\newcommand{\mchrob}{\widehat{\mathcal{H}}}

\newcommand{\fN}{f^\lambda_{N}}

\newcommand{\mcx}{\mathcal{X}}
\newcommand{\blambda}{\bm{\lambda}}

\newcommand{\R}{\mathbb{R}}
\newcommand{\N}{\mathbb{N}}
\newcommand{\E}{\mathbb{E}}
\newcommand{\linspan}{\textnormal{span}}
\newcommand{\sig}{\textnormal{Sig}}
\newcommand{\sigrob}{\widehat{\textnormal{Sig}}}

\newcommand{\Erel}{\mc{E}_{\textnormal{rel}}}
\newcommand{\Tr}{\textnormal{Tr}}
\newcommand{\HS}{\textnormal{HS}}
\newcommand{\dd}{\mathrm{d}}

\DeclareMathOperator*{\argmin}{arg\, min}

\theoremstyle{definition}
\newtheorem{theorem}{Theorem}
\newtheorem{assumption}[theorem]{Assumption}

\theoremstyle{remark}

\newtheorem{definition}[theorem]{Definition}
\newtheorem{remark}[theorem]{Remark}
\theoremstyle{theorem}
\newtheorem{lemma}[theorem]{Lemma}
\newtheorem*{theorem*}{Theorem}
\newtheorem*{definition*}{Definition}
\newtheorem*{fact*}{Fact}
\newtheorem*{corollary*}{Corollary}
\newtheorem{corollary}[theorem]{Corollary}
\numberwithin{equation}{section}

\usepackage{mathtools}
\DeclarePairedDelimiter\abs{\lvert}{\rvert}
\DeclarePairedDelimiter\norm{\lVert}{\rVert}

\usepackage{listings}
\def\mb{\ensuremath\mathbb}

\usepackage[T1]{fontenc}
\usepackage[utf8]{inputenc}
\usepackage[hidelinks]{hyperref}

\usepackage{amsmath}
\usepackage[english]{babel}
\usepackage{mathtools}

\usepackage{listings}
\def\mb{\ensuremath\mathbb}
\usepackage{url}

\begin{document}

\maketitle

\begin{abstract}
    Path signatures are a standard tool for learning functionals on path space. 
    The main underpinnings of their application are universal approximation results. For a wide range of function classes on path space, these results state that any function in such a class can be approximated arbitrarily well by a linear functional of the path signature. 
    Hence, linear functionals of the signature are in principle expressive enough for learning functions on path space. Nonetheless, these results do not provide insights into the statistical question of how much data is needed to identify a suitable signature functional. In this work, we study this question 
    in the framework of kernel ridge regression with the signature kernel. 
    The target is the conditional expectation $\E[Y|X]$, where $X$ is the observed path and $Y$ is a real-valued response.
    Under different assumptions on the moments of the signature kernel and the response, we provide finite-sample and consistency guarantees for the method. We supplement our theoretical analysis with numerical experiments that illustrate our results.
\end{abstract}


\newpage
\tableofcontents
\newpage

\section{Introduction}
\label{sec:introduction}
Path signatures, i.e., the collection $\sig(\bx)$ of all iterated integrals
\begin{equation}
  \label{eq:iterated-integrals}
  \sig^w(\bx) \coloneqq \int_{0 < t_1 < \cdots < t_n < t} \dd \bx^{i_1}_{t_1} \cdots \dd \bx^{i_n}_{t_n}, \quad w = (i_1, \ldots, i_n) \in \{1, \ldots, d\}^n, \ n \ge 1,
\end{equation}
of a path $\bx:[0,t] \to \R^d$ (assumed to be regular enough for the integrals to make sense), have become a standard tool for mathematical and statistical analysis on path space.
We refer to \cite{chevyrev2025primer} for the classical introduction and \cite{mcleod2025signature} for a recent overview article, as well as section~\ref{sec:sign-kern-meth} in this text.
Since Graham's seminal application paper \cite{graham2013sparse} on Chinese handwritten character recognition, signatures have found numerous applications from medicine (see \cite{morrill2020utilization}), weather forecasting (see \cite{dodson2025signature}), quantitative finance (see \cite{cuchiero2025joint}), to dynamic topology and scientific data (see \cite{giusti2023signatures}), among numerous other applications.

Beyond statistical learning, signatures also play an important role in solving stochastic optimal control problems, specifically in the absence of the Markov property.
Indeed, owing to the flow of information in stochastic optimal control, repeated calculations of conditional expectations are a key computational step of many methods for solving such problems, see, e.g., the celebrated Longstaff--Schwartz algorithm for optimal stopping \cite{longstaff2001valuing}.
In terms of statistical learning, this corresponds to \emph{regression}.
We refer to \cite{bank2025stochastic} for a fundamental approach to general stochastic optimal control problems, \cite{kalsi2020optimal,cuchiero2025signature} to applications in optimal execution and portfolio optimization, respectively, and to \cite{BayerPelizzariSchoenmakers2025} for the construction of primal-dual methods in the context of optimal stopping, or, in financial terms, pricing of American options.
In all these cases, the underlying dynamics is non-Markovian or completely data-driven.

Methodologically, signatures are used in several different ways.
At the core lies the observation of universality of linear functionals of the signature with the continuous functions on path-space, see, e.g., \cite{signature_universality_classic} -- in line with the very useful analogy of signatures (or, rather, linear functionals of signatures) as \emph{polynomials on path space}.
Hence, we can use linear functionals of the (truncated) signature for regression on path space.
Additionally, iterated integrals -- i.e., individual components of the signature -- can also be used as features in a machine learning pipeline, e.g., as inputs to neural networks, transformers, random forests, \ldots.
If linearity and associated features (such as the shuffle identity) cannot be taken advantage of, the signature could also be replaced by the log-signature, leading to lower-dimensional features without information-loss.
For more context, we refer to \cite{liao2019learning,bonnier2019deep}, and \cite{bayer2024pricing} for comparisons in the context of American option pricing.

An important development was the introduction of the \emph{signature kernel} \cite{oberhauser_signature_kernel}, linking signature methods with classical kernel methods in machine learning, see, e.g., \cite{scholkopf2001learning}.
In its most basic form, given two paths $\bx:[0,t] \to \R^d$, $\bx': [0,s] \to \R^d$, the signature kernel is defined as
\begin{equation}
  \label{eq:signature-kernel-basic}
  k_{\sig}(\bx,\bx') \coloneqq \sum_{n=0}^\infty \sum_{w \in \{1, \ldots, d\}^n} \sig^w(\bx) \sig^w(\bx')
\end{equation}
-- or its truncated version.
We again refer to section~\ref{sec:sign-kern-meth} for more details on the definition and properties of the signature kernel.
As an alternative to the computationally demanding formula~\eqref{eq:signature-kernel-basic}, \cite{goursat_PDE} provide a Goursat PDE (in the variables $(t,s)$, i.e., for fixed paths $\bx$, $\bx'$) for the signature kernel, thereby allowing to take advantage of the celebrated \emph{kernel trick}, i.e., the observation that kernel methods typically do not require explicit evaluations of the infinite-dimensional feature vector (here: the signature $\sig(\bx)$), but only of the (real-valued)  kernel.
Further developments in the context of the signature kernel include a study of a notion of random Fourier features in \cite{toth2025random}.
Typically, Gaussian processes provide a complementary approach to kernel methods, see \cite{kanagawa2018gaussian} for an overview.
This connection has also been explored in the infinite-dimensional path-space context in \cite{toth2020bayesian}.

We are interested in the \emph{regression problem} on path-space.
Let $\mathcal{X}$ denote a space of paths, i.e., functions $[0,t] \to \R^d$, $t > 0$.
To make notions of integration precise, we will work in a \emph{rough path setting}, see Section~\ref{sec:sign-rough-paths} for details.
In fact, we will choose $\mathcal{X}$ to be a Polish space.
We consider a path-valued random variable $X$ together with a real random variable $Y$, and want to compute
\begin{equation}
  \label{eq:regression-problem}
  f^{\ast}(\bx) \coloneqq \mathbb{E}[Y \mid X = \bx], \quad \bx \in \mathcal{X}.
\end{equation}
It is well understood that signatures allow us to solve the regression problem in the sense that functions on $\mathcal{X}$ can be approximated arbitrarily well by linear functionals on the signature in the $L^p$-sense.
This was established in different settings by \cite{BayerPelizzariSchoenmakers2025} (building on the robust universal approximation result of \cite{chevyrev_oberhauser_robust}) as well as \cite{ceylan2025globaluniversalapproximationbrownian} (building on the universal approximation result due to \cite{Cuchiero2026} w.r.t.~weighted norms).

\emph{Expressivity}, however, is only the starting point.
In practice, the regression problem needs to be solved based on finitely many \emph{samples} $(X_i, Y_i)_{i=1}^N \sim (X,Y)$, and statistical analysis aims to quantify the error, i.e., the difference between $f^{\ast}$ and its estimator $f_N^\lambda$ -- where $N$ denotes the number of samples and $\lambda>0$ denotes a regularization parameter, see Section~\ref{sec:sign-kern-meth} for details.
Note that $f_N^\lambda$ is itself random, as it depends on the \emph{training data} $(X_i, Y_i)_{i=1}^N$.

The seminal paper \cite{fermanian2022functional} studies the problem under the assumption that the unknown function $f^{\ast}$ is given as a finite linear combination of iterated integrals.
Left with estimating the truncation degree and the coefficients, an error estimate for the squared prediction error of order $C N^{-1/2}$ is derived.
In the related work \cite{bleistein2023learning}, this method is taken up and applied to the problem of learning the solution map of stochastic differential equations (SDEs), again as a linear combination of signature terms.
Again, error estimates in terms of squared prediction errors are provided.

Finite-sample error estimates with high probability in terms of prediction errors are provided in \cite{horvath2026fast} in the case of OLS (as above), but also for LASSO regression and logistic regression with signatures.
Logistic regression with signatures including adaptively chosen truncation is also developed and analyzed in \cite{zeng2025adaptively}.
Similarly, \cite{bayer2026local} provide finite-sample error estimates with high probability for Nadaraya--Watson type local regression on path space based on rough path metrics and signatures. 

In order to robustify signature methods, \cite{chevyrev_oberhauser_robust} suggest to replace the signature by a normalized version, referred to as the \emph{robust} signature.
Different to the standard signature kernel, the corresponding robust signature kernel is bounded.
The problem of regression with stochastic processes using robust signatures is studied in \cite{schell_signatures}. They provide consistency guarantees for kernel ridge regression with the robust signature kernel but do not provide rates.

In this work, we provide finite-sample theory for signature kernel ridge regression estimators of $f^{\ast}$.
More specifically, let $\mathcal{H}$ denote the RKHS associated to the signature kernel $k_{\sig}$.
Given i.i.d.~samples $(X_i, Y_i)_{i=1}^N$, define
\begin{equation}
  \label{eq:estimator}
  f^\lambda_N\coloneqq \argmin_{f \in \mathcal{H}} \frac{1}{N} \sum_{1\leq i \leq N} (Y_i-f(X_i))^2 +\lambda \|f\|_{\mch}^2,
\end{equation}
see Section \ref{sec:kern-ridge-regr} for more details.
We measure the error in terms of \emph{excess risk}
\begin{equation*}
  \mc{E}_{\textnormal{rel}} [f] \coloneqq \int_{\mcx\times \R} \abs{f(\bx) - y}^2 \nu(\dd \bx, \dd y) - \int_{\mcx\times \R} \abs{f^{\ast}(\bx) - y}^2 \nu(\dd \bx, \dd y),
\end{equation*}
where $\nu$ denotes the joint distribution of $X$ and $Y$.
In particular, when considering $\mc{E}_{\textnormal{rel}} [f^\lambda_N]$, the randomness of $f^\lambda_N$ -- which depends on the \emph{training data} $(X_i, Y_i)_{i=1}^N$ -- persists in the random excess risk.
Our analysis follows the classical approach by \cite{Caponnetto2007OptimalRates}, with different adaptations required for different cases.
Assuming that $f^{\ast} \in \mathcal{H}$, we obtain different estimates under different scenarios.

\paragraph{Light tailed case}

Under strong moment assumptions on $k_{\sig}(X,X)$ as well as strong conditional moment assumptions on $Y - f^{\ast}(X)$ (see Assumption~\ref{assumption:finitesample:bernstein}), we can use Bernstein's inequality as the basis of the finite-sample error bound. In this case, we obtain for $\lambda = \lambda_N \sim 1/\sqrt{N}$ with probability at least $1-\eta$ for large enough $N$,
\begin{equation*}
  \mc{E}_{\textnormal{rel}} [f^{\lambda_N}_N] \le C \sqrt{\frac{\log(6/\eta)}{N}},
\end{equation*}
with a constant $C$ depending on $\norm{f^{\ast}}_{\mathcal{H}}$, see Theorem~\ref{theorem:rates:Bernstein}.

\paragraph{Heavy tailed case}

It turns out that the moment conditions required for Bernstein's inequality are not satisfied for the signature kernel in most cases of practical interest.
Specifically, it turns out, see Lemma~\ref{lemma:BM-kernel-moments} together with~\eqref{eq:lb-moments-strat-bm}, that the signature kernel does not satisfy the conditions of Assumption~\ref{assumption:finitesample:bernstein} when the driving path $X$ is Brownian motion.
However, we can replace Bernstein's inequality by the Fuk-Nagaev inequality, under the much weaker Assumption~\ref{assumption:finitesample:heavytail}.
Theorem~\ref{theorem:rates:heavytail} then provides -- again for $\lambda_N \sim 1/\sqrt{N}$ and large enough $N$--
\begin{equation*}
  \mc{E}_{\textnormal{rel}} [f^{\lambda_N}_N] \le C \sqrt{\frac{\log(6c/\eta)}{N}},
\end{equation*}
with different constants $C,c>0$ and a stronger lower bound condition on $N$. 
\paragraph{Robust signature}

Since the robust signature kernel is bounded, we can analyze robust signature kernel regression using the error estimates of \cite{Caponnetto2007OptimalRates}, yielding
a $1/\sqrt{N}$ rate on the excess risk, as in the light and heavy tailed cases. A key difference to the other two cases is though that no distribution dependent moment assumptions need to be imposed on the robust signature kernel, as it is bounded by construction. 

For all three cases we also provide consistency guarantees that do not rely on $f^{\ast}\in \mch$.

The theoretical analysis is accompanied by numerical examples (see Section~\ref{sec:numer-exper}), showcasing the asymptotic behavior of the excess risk as well as the optimal penalization coefficient $\lambda$ in three examples.
In Section~\ref{sec:one-dimensional-sde}, $f^{\ast}$ denotes the solution map of an SDE -- mapping the driving Brownian path to the solution at terminal time.
We then recover a signature-based model of a financial market in Section~\ref{sec:hest-sign-model}.
Finally, in Section~\ref{sec:tser}, we apply signature kernel ridge regression to a database \cite{TSERTan2020} of regression problems with time series data.

\section{Signatures and kernel methods}
\label{sec:sign-kern-meth}
In this section, we introduce the problem setup. In Section \ref{sec:sign-rough-paths}, we rigorously introduce the concepts of signatures and rough paths, followed by an introduction to (signature) kernel ridge regression in Section \ref{sec:kern-ridge-regr}. We close with a discussion of the signature RKHS in Section \ref{sec:signature-rkhs}.
\subsection{Signatures and rough paths}
\label{sec:sign-rough-paths}
For $d\in\mb{N}$, we denote the standard basis of $\R^d$ by $\{e_1,\dots,e_d\}$. We define further the tensor algebra by 
\begin{align}\notag
    T(\mb{R}^d)\coloneqq \left\{v=(v_0,v_1,\dots,v_k,\dots)\in \prod_{k\in \mb{N}} (\mb{R}^d)^{\otimes k}:\exists K\in \mb{N} \textnormal{   s.t.  } v_k=0
    \textnormal{  for all }k\geq K \right\},
\end{align}
and the extended tensor algebra by $T((\mb{R}^d))\coloneqq 
    \prod_{k\in \mb{N}} (\R^d)^{\otimes k}$,
using the convention $(\mb{R}^d)^{\otimes 0}=\mb{R}$.
We also write $T_1((\mb{R}^d))\coloneqq\{v\in T((\mb{R}^d)):v_0=1\}$.
Define for $v=(v_0,v_1,\dots),v'=(v_0',v_1',\dots)\in T((\R^d))$ the formal inner product
\begin{align}
    \langle v, v' \rangle_{T^2((\R^d))} \coloneqq \sum_{k\in \N}
    \langle v_k, v'_k \rangle_{(\R^d)^{\otimes k }}.
\end{align}
We then set 
\begin{align}
    T^2((\R^d)) \coloneqq 
    \left \{ v\in T((\R^d)):\,
    \|v\|_{T^2((\R^d))}<\infty
    \right\}.
\end{align}
Further, for $\phi:\N\to \R_+ \setminus \{0\}$ and $v,v'\in T((\R^d))$, we set
\begin{align}
    v^\phi \coloneqq (\sqrt{\phi(0)}v_0,\sqrt{\phi(1)}v_1,\dots)\in T((\R^d))
\end{align} 
and $\langle v, v' \rangle_{T^{2,\phi}((\R^d))}\coloneqq \langle v^\phi, v'^\phi \rangle_{T^2((\R^d))}$. Then, we define $T^{2,\phi}((\R^d))$ analogously by
\begin{align}
    T^{2,\phi}((\R^d)) \coloneqq 
    \left \{ v\in T((\R^d)):\,
    \|v\|_{T^{2,\phi}((\R^d))}<\infty
    \right\}.
\end{align}
Further details on the truncated tensor algebra can be found in \cite[Chapter 7.2.1]{Friz_Victoir_2010}.
\paragraph{Geometric rough paths.} 
In this work, we assume all paths to start at zero.
We denote for $\alpha \in (\frac{1}{3},1]$ the space of $\alpha$-Hölder continuous paths by $C^\alpha ([0,T],\mb{R}^d)$, i.e., $X\in C^\alpha ([0,T],\mb{R}^d)$ if 
\begin{align}\notag
    \|X\|_{\alpha,[0,T]}\coloneqq \sup_{0\leq s \leq t\leq T} \frac{\|X_t-X_s\|}{|t-s|^{\alpha}}
    <\infty.
\end{align}
Here and in the following, we denote by $\|\cdot \|$ the Euclidean norm on $\mb{R}^d$. 
If there is no ambiguity about the time scale, we make the dependence implicit and write $\|X\|_{\alpha}=\|X\|_{\alpha,[0,T]}$. Further, we denote the $[0,T]$-simplex by
\begin{align}
    \Delta_{[0,T]}^2 \coloneqq \left\{(s,t)\in [0,T]^2:\,0\leq s \leq t \leq T
    \right\}.
\end{align}
Given a two parameter function $\mathbf{X}:\Delta_{[0,T]}^2 \to T^{\leq 2}(\mb{R}^d)\coloneqq \R \times (\R^d)\times (\R^d)^{\otimes 2}$, so $\mathbf{X}_{s,t}=(1,\mathbf{X}^{(1)}_{s,t},\mathbf{X}^{(2)}_{s,t})$, we define the rough path norm $\vert\!\vert\!\vert \cdot \vert\!\vert\!\vert_{\alpha}$ by
\begin{align}
    \vert\!\vert\!\vert \mathbf{X} \vert\!\vert\!\vert_{\alpha}
    \coloneqq 
    \underset{\ell = 1,2}{\textnormal{max}}
    \,
    \left(
    \sup_{0\leq s \leq t\leq T} \frac{\|\mathbf{X}_{s,t}^{(\ell)}\|}{|t-s|^{\ell \alpha}}
    \right)^{1/\ell}.
\end{align}
The space of geometric $\alpha$-Hölder rough paths $\mathscr{C}^{\alpha}_g([0,T],\mb{R}^{d})$
is then defined as the limit of the canonical rough path lifts of Lipschitz-continuous paths $X:[0,T]\to\mb{R}^d$ w.r.t. the rough paths norm $\vert\!\vert\!\vert \cdot \vert\!\vert\!\vert_{\alpha}$. For further discussion we refer the interested reader to \cite[Chapter 9]{Friz_Victoir_2010}.
\begin{remark}
    Note that in our framework the rough path sees only the increments, which is natural given that we assume all paths to start in zero. One can also consider rough paths with an initial condition, see \cite{Friz_Victoir_2010}.
\end{remark}

\paragraph{Signatures of geometric rough paths.} Given a path $\mathbf{X}\in \mathscr{C}^{\alpha}_g([0,T],\mb{R}^{d})$ for $\alpha \in (1/3,1]$, we define the rough-path signature path as the unique path $\mathbf{X}^{<\infty}:\Delta^2_{[0,T]}\to T((\mb{R}^d))$
from Lyons' extension Theorem taking values in the group of group-like elements $G((\R^d))\subset T((\R^d))$, see \cite{Lyons1998}, Theorem 3.7. 

In the following, we refer to the signature of a rough path as $\sig(\bX)\coloneqq \bX^{<\infty}_{0,T}$.
The path is uniquely characterized up to tree-like equivalences
by its signature together with its initial value, see \cite{Boedihardjo2016}.
When considering $\alpha$-geometric rough paths augmented by a time component, the signature of the geometric rough path lift of the time augmented process together with its initial value uniquely identifies the path among the class. In the following, we will denote the space of time augmented $\alpha$-Hölder geometric rough paths by $\widehat{\mathscr{C}}^\alpha_g([0,T],\R^{d+1})$, i.e., the space of lifts of $\widehat{X}_t=(t,X_t)$, where $X$ is $\alpha$-Hölder continuous.

The signature is indexed by multi-indices, i.e., elements of $\cup_{n\in \N}[d]^n$. Those multi-indices are referred to as words. We write $W\coloneqq \cup_{n\in \N}[d]^n$ for the set of words and $W_n\coloneqq [d]^n$ for the set of words of length $n$.

\paragraph{Stopped geometric rough paths.} 
The concept of rough paths can be extended to the setting where the time-scale of paths is not fixed, to which one refers to as stopped rough paths \cite{BayerHagerRiedelSchoenmakers2023, BayerPelizzariSchoenmakers2025}. 
We now give the definition of stopped geometric rough paths.
\begin{definition}[{\cite[Definition 2.1]{BayerPelizzariSchoenmakers2025}}]\label{def:stoppedRP} 
For any $\alpha \in (\frac{1}{3},1]$ and $T>0$, the space of stopped $\widehat{\mathscr{C}}^{\alpha}_g-$paths on $[0,T]$ is defined by the disjoint union 
    $$\Lambda^{\alpha}_T:= \bigcup_{t\in [0,T]}\widehat{\mathscr{C}}^{\alpha}_g([0,t],\R^{d+1}).$$ Moreover, we equip the Polish space $\Lambda_T^{\alpha}$ with the final topology induced by the map $$\phi^\ast:[0,T] \times \widehat{\mathscr{C}}^{\alpha}_g([0,T],\R^{d+1}) \longrightarrow \Lambda_T^{\alpha},\quad \phi^\ast(t,\mathbf{x}) = \mathbf{x}|_{[0,t]}.$$
\end{definition} 
As observed in \cite{BayerPelizzariSchoenmakers2025}, the space $\Lambda_T^\alpha$ can be equipped with a metric $d_\Lambda$ whose induced topology is the one of Definition \ref{def:stoppedRP}.

\subsection{Kernel ridge regression}
\label{sec:kern-ridge-regr}

Consider two random variables $X\in \mc{X}$ and $Y\in \R$ with joint distribution $\nu\in \mc{P}(\mc{X}\times \R)$ for $\mc{X}$ a Polish space endowed with its Borel sigma-algebra. We denote the first and second marginals of $\nu$ by $\nu_{\mcx},\nu_{\R}$.
Given i.i.d. realizations $(x_i,y_i)_{1\leq i \leq N}$ of the random variables $X$ and $Y$, the problem of regressing $Y$ on $X$ consists in using the observations to identify a function $f:\mc{X}\to \R$ that is a good estimator of the conditional expectation $\mb{E}[Y|X]$, i.e.,
\begin{align}\label{eq:approxcond}
    f(x)\approx \mb{E}[Y|X=x] \eqqcolon f^{\ast}(x).
\end{align}
Given some Hilbert space $H$ of functions $f:\mcx \to \R$, in ridge regression one does so by minimizing across $f\in H$ the $L^2$ regularized empirical squared loss functional
\begin{align}\label{eq:empiricallossridge}
    f\mapsto \frac{1}{N} \sum_{1\leq i \leq N} (y_i-f(x_i))^2
    +\lambda \|f\|_H^2,
\end{align}
where $\lambda>0$. This approach is also known as Tikhonov regularization \cite{Tikhonov1963}. 
Observe that the problem of minimizing 
\eqref{eq:empiricallossridge} across $H$ admits a unique solution provided the mapping \eqref{eq:empiricallossridge} is lower semicontinuous, as it 
is strongly convex.
When dropping the regularization, i.e., for $\lambda=0$, the problem of minimizing the finite-sample objective \eqref{eq:empiricallossridge} 
across $H$ would still admit a solution, though uniqueness and stability are not guaranteed.

Typically, one considers not a generic Hilbert space, but a \textit{reproducing kernel Hilbert space} (RKHS), which we denote in the following by $\mch$. 
Those spaces are characterized by having continuous evaluation functionals. That is, for any fixed $x\in \mc{X}$ the map
\begin{align}\label{eq:reproducing}
    k_x:\mch\to \R, \quad k_x f \coloneqq f(x)
\end{align}
is continuous w.r.t. the $\mch$-norm. In particular, then also \eqref{eq:empiricallossridge} is lower semicontinuous.

Note that by Riesz, we have $k_x f = \langle k_x,f\rangle_\mch$, where by an abuse of 
notation we denote both the continuous functional $k_x:\mch\to \R$ as well as its 
Riesz representative by $k_x$. In particular, one can define a kernel associated to 
the Hilbert space $\mch$ for $x,x'\in \mcx$ by
\begin{align}
    k(x,x') \coloneqq \langle k_x,k_{x'}\rangle_{\mch} .
\end{align}
Observe that the kernel can be interpreted both as a map $k(\cdot,\cdot): \mcx \times \mcx \to \R$ as well as a collection of maps $\{k_x\}_{x\in \mcx}$, see \eqref{eq:reproducing}. Due to \eqref{eq:reproducing}, one says that the kernel \textit{reproduces} the Hilbert space. In the following, we will also use the abbreviation $k_x= k(x,\cdot)$.

As was observed initially in \cite{Aronszajn1950}, see also \cite[Theorem 2.5]{kernel_mean_embeddings_review}, there is a one-to-one correspondence between 
reproducing kernel Hilbert spaces and kernels.
\begin{theorem}[{\cite[Theorem 2.5]{kernel_mean_embeddings_review}}]\label{theorem:RKHS-kernel}
    Let $\mc{X}$ be a nonempty set. For every positive definite function $k(\cdot,\cdot):\mc{X} \times \mc{X} \to \R$, there exists a unique (up to isometric isomorphism) reproducing kernel Hilbert space (RKHS) with $k$ as its reproducing kernel. Conversely, the reproducing kernel of an RKHS is unique and positive definite.
\end{theorem}

In the context of ridge regression, the interest in the kernel stems from the fact that the minimizer
\begin{align}\label{eq:flamdbdaN}
    f^\lambda_N\coloneqq \argmin_{f\in \mch} \frac{1}{N} \sum_{1\leq i \leq N} (y_i-f(x_i))^2
    +\lambda \|f\|_{\mch}^2
\end{align}
is an element of the finite dimensional subspace $\linspan \left(\{k_{x_i}\}_{1\leq i \leq N}\right)\subset \mch$. 
This classical result is known as the representer Theorem, and holds also in more general frameworks, see \cite[Theorem 1]{Scholkopf2001GeneralizedRepresenter}.

Hence, the minimizer $f^\lambda_N$ can be identified by minimizing \eqref{eq:empiricallossridge} across the finite dimensional subspace $\linspan\left( \{k_{x_i}(\cdot)\}_{1\leq i \leq N}\right)$, which is a fully tractable minimization problem. Invoking the reproducing property, we find that the coefficient $\alpha^\ast$ associated to the representation $f^\lambda_N = \sum_{1\leq i \leq N} \alpha^\ast_i k_{x_i}$ is given by
\begin{align}\label{eq_alpha_ridge_Rd}
    \alpha^\ast
    =
    (K+\lambda N )^{-1} \,\hat{y}
    \in \mb{R}^N,
\end{align}
where $K=\left(k(x_i,x_j)\right)_{1\leq i,j\leq N} \in \R^{N\times N}$ and $\hat{y}=(y_i)_{1\leq i \leq N}\in \mb{R}^N$. Therefore, to evaluate $f^\lambda_N$, we have to compute only kernels, i.e. quantities of the form $k(x,x')$. Because of this, in our setting one refers to ridge regression as \textit{kernel} ridge regression. 

We want to use kernel ridge regression for learning functions of (rough) paths using 
the \textit{signature kernel} \cite{oberhauser_signature_kernel, WeightedSignatureKernels}, whose RKHS is highly expressive, see Section 
\ref{sec:signature-rkhs}. Let $\phi: \mb{N}\to\mb{R}_+\setminus \{0\}$ fulfill for all $C>0$ that
\begin{align}\label{eq_condition_phi}
        \sum_{k\geq 0} \frac{C^k\phi(k)}{((k/3)!)^2}<\infty.
\end{align}
Following \cite{oberhauser_signature_kernel, WeightedSignatureKernels}, given  two $\alpha$-Hölder continuous geometric (stopped) rough paths $\bx,\bx'\in \Lambda_T^\alpha$, $\alpha\in (\frac{1}{3},1]$, the $\phi$-signature kernel evaluated at $\bx,\bx'$ is defined by
\begin{equation}\label{eq:sigkernel}
    k_{\sig,\phi}(\bx,\bx') \coloneqq 
        \langle \sig(\bx),\sig(\bx') \rangle_{T^{2,\phi}((\R^{d+1}))}
        = \sum_{k\geq 0} \phi(k) \langle \sig(\bx)^{(k)},\sig(\bx')^{(k)}\rangle_{(\R^{d+1})^{\otimes k}},    
\end{equation}
which is well defined due to \eqref{eq_condition_phi}, see, e.g., 
\cite[Exercise 4.6]{friz_hairer}
and \cite[Proposition 2.9]{WeightedSignatureKernels} for a proof in the bounded variation setting.

Note also that in this setup, while enjoying the expressivity of the class of linear functionals of the signature, we do not necessarily have to compute signatures explicitly to compute $\fN$. Indeed, we only have to compute signature kernels. As 
the dimension of the signature level scales exponentially, in high 
dimensions, this is computationally lighter due to the existence of exact efficient PDE-based algorithms for the case of constant $\phi$, see \cite{goursat_PDE}. This is referred to as a \textit{kernel trick}, as one can avoid the costly computation of the signature by efficiently computing the signature kernel.

While for non-constant $\phi$ no kernel trick is directly available, one can still employ a Monte Carlo scheme to reduce the problem of computing $k_{\sig,\phi}$ to the problem of computing standard signature kernels if $\phi(k)$ corresponds to the moments of a random variable, see \cite{WeightedSignatureKernels}.
\subsection{The signature RKHS}
\label{sec:signature-rkhs}
In kernel ridge regression, the functions that are available to estimate $f^{\ast}$ are the elements of the RKHS associated to the kernel. In the case of the $\phi$-signature kernel, those are linear functionals of the signature. 
\begin{lemma}\label{lemma_RKHS}
    Assume $\phi:\N\to \R_+ \setminus \{0\}$ fulfills \eqref{eq_condition_phi}.
    For $\alpha \in (\frac{1}{3},1]$, $T>0$,
    consider the stopped rough-path space $\mcx=\Lambda_T^\alpha$ or the time augmented rough path space $\mcx =
    \widehat{\mathscr{C}}^{\alpha}_g([0,T],\R^{d+1})$. Then, the RKHS $\mch_\phi$ associated to the $\phi$-signature kernel $k=k_{\sig,\phi}$ is given as the closure of the function space
    \begin{equation}
        \widetilde{\mch}_{\phi} \coloneqq \textnormal{span} \{k_{\bx} :\bx'\mapsto k(\bx,\bx'): \quad \bx \in \mcx\},
    \end{equation}
    equipped with the scalar product $\langle \cdot, \cdot \rangle_{\mch_\phi}$, defined for $k_{\bx},k_{\bx'}$ by
    \begin{align}
        \langle k_\bx, k_{\bx'}\rangle_{\mch_\phi} \coloneqq k(\bx,\bx')
    \end{align}
    and extended to $\widetilde{\mch}_{\phi}$ by linearity.
    
    Further, we have 
    \begin{equation*}
        \mch_\phi= \big\{ f:\mcx \to \R: \quad \,f(\bx) = \langle \ell_f,\sig(\bx)\rangle_{T^{2,\phi}((\R^{d+1}))} 
        ,\, \ell_f\in T^{2,\phi}((\R^{d+1}))
        \big\},
    \end{equation*}
    
    and, given $f\in \mch_\phi$, that
    \begin{equation}\label{eq:norm:f}
        \|f\|_{\mch_\phi} = \inf 
        \left\{\|\ell_f\|_{T^{2,\phi}((\R^{d+1}))}:\quad 
        f(\cdot)=\langle \ell_f,\sig (\cdot) \rangle_{T^{2,\phi}((\R^{d+1}))}, \,
        \ell_f \in T^{2,\phi}((\R^{d+1}))
        \right\}.
    \end{equation}
\end{lemma}
The result and its proof can be found in \cite[Theorem 4.21]{steinwart2008support}, where the setting is more general.
\begin{remark}
    Note that, since for $\mcx=\Lambda_T^\alpha$ one can show that there exists for each
    $f\in \mch_\phi$ a unique representative $\ell_f \in T^{2,\phi}((\R^{d+1}))$, the infimum in \eqref{eq:norm:f} can be dropped. This is not the case for $\mcx =\widehat{\mathscr{C}}^{\alpha}_g([0,T],\R^{d+1})$, as the time-channel is path independent.
\end{remark}

Observe that given the construction, it is clear a priori that elements of $\mch_\phi$ are linear functions of the signature. To understand heuristically why also for every $\ell_f \in T^{2,\phi}((\R^{d+1}))$ we have $\langle \ell_f,\sig (\cdot) \rangle_{T^{2,\phi}((\R^{d+1}))}\in \mch$, observe that one can decompose every $\ell_f \in T^{2,\phi}((\R^{d+1}))$ as $\ell_f=\ell_f^\sig + \ell_f^o$, where $\ell_f^\sig \in T^{2,\phi}((\R^{d+1}))$ is in the closure of the span of signatures and $\ell_f^o \in T^{2,\phi}((\R^{d+1}))$ is orthogonal to signatures, i.e.,
\begin{equation*}
    \langle \ell_f^o,\sig(\cdot)\rangle_{T^{2,\phi}((\R^{d+1}))}=0 \quad
    \textnormal{on }\mcx.
\end{equation*}
Hence, $\langle \ell_f,\sig (\cdot) \rangle_{T^{2,\phi}((\R^{d+1}))}
= \langle \ell_f^\sig,\sig (\cdot) \rangle_{T^{2,\phi}((\R^{d+1}))}\in \mch$.

The significance of this characterization for $\mch_\phi$ lies in the fact
that linear functionals of the signature are universal approximators for $L^p$ functionals, in particular $\mch_\phi$ is dense in $L^p$ spaces. This is the content of the following
result due to Ceylan and Prömel, building on work of Bayer et al. \cite{BayerPelizzariSchoenmakers2025}.
\begin{theorem}[{$L^p$-universal approximation Theorem on $\Lambda_T^\alpha$, \cite[Theorem 3.13]{ceylan2025globaluniversalapproximationbrownian}}]
\label{theorem_approx_global}
Define $\psi:\Lambda_T^\alpha \to \mb{R}$ by
\[
\psi(\bx_{[0,t]})
:=
\exp\!\left(
\beta \, \| \bx_{[0,t]} \|^{\gamma}_{\textnormal{cc},\alpha}
\right),
\qquad
\bx_{[0,t]} \in \Lambda_T^\alpha
\]
for some $\beta>0$, $\gamma \geq \lfloor 1/\alpha \rfloor$, where $\| \cdot \|_{\textnormal{cc},\alpha}$ denotes the Carnot-Carathéodory rough path norm. Consider a probability measure $\mu \in \mc{P}(\Lambda_T^\alpha)$.
Let $p>1$, and assume that
\begin{align}\label{eq_integrability_approx}
\int_{\Lambda_T^\alpha} \psi^p (\bx) \, \mu(\mathrm{d}\bx) < \infty .
\end{align}
Then, for every $f \in L^p(\Lambda_T^\alpha,\mu)$ and for every $\varepsilon>0$, there exists a functional $f_\ell =\langle \ell_f, \sig (\cdot)\rangle$ with $\ell_f\in T(\R^{d+1})$, such that
\[
\|f - f_\ell\|_{L^p(\Lambda_T^\alpha,\mu)} < \varepsilon.
\]
\end{theorem}
The proof of Theorem \ref{theorem_approx_global} relies on a modification of a novel weighted version of the Stone–Weierstrass Theorem, see \cite[Theorem 5.4]{Cuchiero2026}.
Importantly, as checked in \cite[Section 4]{ceylan2025globaluniversalapproximationbrownian}, time augmented Brownian motion fulfills the integrability condition \eqref{eq_integrability_approx} for suitable choices of the parameters, leading to the following result.
\begin{theorem}[{\cite[Corollary 4.3]{ceylan2025globaluniversalapproximationbrownian}}]\label{theorem_approx_B_sig}
  Let $\alpha \in (1/3,1/2)$, let $W$ be a Brownian motion, $\widehat W=(\cdot,W)$ be the time-extended Brownian motion and $\widehat{\mb{W}}$ be the corresponding time-extended Brownian geometric rough path. Denote by $\mathbb H^p$ the space of $(\mathcal F^W_t)_{t\geq 0}$-progressively measurable processes $A$ fulfilling
  \begin{equation*}
      \|A\|_{\mathbb H^p}^p:=\E\Bigl[\int_0^T|A_t|^p\, dt\Bigr]<\infty.
    \end{equation*}
    Then, the following two statements hold.
  \begin{enumerate}
    \item[(i)] Let $f(\widehat{\mathbb{W}})\in L^p(\Omega;\mb{P})$ with $f\colon \widehat{\mathscr{C}}^{\alpha}_g([0,T],\R^{d+1}) \to \R$. Then, for every $\varepsilon>0$ there exists a linear function $f_\ell =\langle \ell_f, \sig (\cdot)\rangle$ with $\ell_f\in T(\R^{d+1})$, such that
    \begin{equation*}
      \E[|f(\widehat{\mathbb{W}})-f_\ell(\widehat{\mb{W}})|^p]<\varepsilon.
    \end{equation*}
    \item[(ii)] Let $f(\widehat{\mathbb{W}}_{[0,\cdot]})\in\mathbb H^p$ with $f\colon \Lambda_T^\alpha \to \R$. Then, for every $\varepsilon>0$ there exists a linear function $f_\ell =\langle \ell_f, \sig (\cdot)\rangle$ with $\ell_f\in T(\R^{d+1})$, such that
    \begin{equation*}
      \E\Bigl[\int_0^T|f(\widehat{\mathbb{W}}_{[0,t]})-f_\ell(\widehat{\mb{W}}_{[0,t]})|^p\,d t\Bigr]<\varepsilon.
    \end{equation*}
  \end{enumerate}
\end{theorem}

\section{Theoretical guarantees for approximation and finite-sample error}\label{sec:theory}
We now provide a theoretical analysis for signature-kernel ridge regression. The proofs of the results stated in this section can be found in the appendix, see Sections \ref{sec:approx-error:proofs}, \ref{section_sub_finite_sample_error}, and \ref{sec:heavytailed:proofs}. For clarity of exposition, we restrict ourselves to the standard case $\phi=1$, though our results also apply to the weighted case. 
In the following, we will always assume $\mc{X}=\Lambda_T^\alpha$
or $\mc{X}=\widehat{\mathscr{C}}^{\alpha}_g([0,T],\R^{d+1})$
and denote the RKHS associated to the signature kernel by $\mc{H}$.
Throughout this work we will make the following integrability and density assumption. 
\begin{assumption}\label{ass_kernel_diagonal}
    The measure $\nu\in \mc{P}(\mcx\times \R)$ is such that $\nu_{\R }\in \mc{P}_2(\R)$ and 
    \begin{align}\label{ass_kernel_diagonal_eq}
        \|k_{\sig}\|^2_{L^2(\nu_{\mc{X}})}\coloneqq 
        \int_{\mcx}k_{\sig}(\bx,\bx)\,\,\nu_{\mc{X}}(\mathrm{d}\bx)<\infty.
    \end{align}
    Further, the RKHS $\mch$ is dense in $L^2(\nu_\mcx)$ and $\textnormal{dim}(L^2(\nu_\mcx))=\infty$.
\end{assumption}
\begin{remark}
    The non-restrictive dimension assumption $\textnormal{dim}(L^2(\nu_\mcx))=\infty$ is made purely for clarity of exposition. All our results hold and are proven in the same way if $\textnormal{dim}(L^2(\nu_\mcx))<\infty$.
\end{remark}
By \cite[Lemma 2.3]{SteinwartScovel2012}, the integrability assumption \eqref{ass_kernel_diagonal_eq} ensures that the embedding $\mc{H} \hookrightarrow L^2(\nu_\mcx)$ is compact.
Observe that Assumption \ref{ass_kernel_diagonal} holds for the law of time augmented $d$-dimensional Brownian motion on $[0,T]$ for $T\in (0,\infty)$.

In the following, for all choices of space $\mc{X}$, we write
\begin{align}\label{eq_ridge_regression}
    f^{\lambda} \coloneqq \underset{f\in \mch}{\argmin} \Big\{ \int_{\mc{X}\times \mb{R}} |y-f(\bx)|^2 \nu(\mathrm{d}\bx,\mathrm{d}y)+\lambda \|f\|^2_{\mch}\Big\}.
\end{align}
Define the operator
\begin{align}\label{eq:defT}
    \mc{T}:\,L^2(\nu_{\mc{X}})\to \mch,
    \quad
    \mc{T} f \coloneqq \int_{\mcx}k_{\bx'}(\cdot) f(\bx')\,\nu_{\mc{X}}(\mathrm{d}\bx').
\end{align}
Since Assumption \ref{ass_kernel_diagonal} ensures that $\mch$ is compactly embedded into $L^2 (\nu_{\mc{X}})$, the operator $\mc{T}$ is compact.
Further, both its restriction to $\mch$ and its embedding into $L^2(\nu_\mcx)$ are positive,  self adjoint, and compact, see \cite{SteinwartScovel2012}.
In \eqref{eq:defT}, integrals are to be understood in a Bochner sense.

Under Assumption \ref{ass_kernel_diagonal}, the minimizer $f^\lambda$ can be represented in terms of $\mc{T}$ by
\begin{align}\label{eq_min_ridge_general}
    f^{\lambda} = (\mc{T}+\lambda)^{-1} \,\mc{T} f^{\ast}.
\end{align}
For a detailed proof of \eqref{eq_min_ridge_general} in a similar setting, see for example \cite[Proposition 3.3]{BayerBelomestnyButkovskyEtAl2024}.

Given $N$ samples $\{z_i\}_{1\leq i \leq N}=\{(\bx_i,y_i)\}_{1\leq i \leq N}\subset \mcx \times \R$, we denote the finite-sample analogue of \eqref{eq_ridge_regression} by
\begin{align}
    f_{N}^\lambda \coloneqq \underset{f\in \mch}{\argmin} \Big\{ 
    \frac{1}{N} \sum_{1\leq i \leq N} |y_i-f(\bx_i)|^2+\lambda \|f\|^2_{\mch}\Big\} \in \mch,
\end{align}
where the dependence on $\{z_i\}_{1\leq i \leq N}$ is made implicit in the notation.

The error analysis for ridge regression consists in quantifying the approximation error and the finite-sample error. The former refers to the error occurring due to $f^{\lambda}$ differing from the conditional expectation $f^{\ast}(\cdot)=\mb{E}[Y|X=\cdot]$. The latter refers to the total error occurring in the estimation of $f^{\ast}$ by means of the stochastic finite-sample estimator $f^{\lambda}_N$. Note that the finite-sample error depends on the approximation error.

The results on the analysis of the approximation error are contained in Section \ref{sec:approx-error}. The finite-sample results obtained under respectively a Bernstein and a Fuk-Nagaev condition on the moments of the signature kernel are the contents of Sections \ref{sec:finite-sample:bernstein} and \ref{sec:finite-sample:fuk-nagaev}.
Those concentration inequalities, i.e., Bernstein and Fuk-Nagaev, are discussed in Section \ref{sec:concentration}.

\subsection{Analysis of the approximation error}
\label{sec:approx-error}

To state the guarantee for the approximation error, we need some preliminaries on $\mc{T}$. 
Due to  $\mc{T}$ being positive, self adjoint, and compact, as a consequence of Lemma \ref{lemma:steinwart-scovel} stated in Section \ref{sec:approx-error:proofs} (essentially the spectral Theorem), we find that there exists an orthonormal basis (ONB) $([e_i]_{\sim})_{i \in \mb{N}}$ of $L^2(\nu_{\mc{X}})$ induced by an orthogonal system $(e_i)_{i \in \mb{N}}$ in $\mc{H}$ with associated (strictly) positive, non-increasing eigenvalues $(\mu_i)_{i\in \mb{N}}$, such that
\begin{align}\label{eq_T_k_L2_classical_1}
    [\mc{T} f]_{\sim} = \sum_{i\in \mb{N}}\mu_i \langle f,[e_i]_{\sim}\rangle_{L^2(\nu_{\mc{X}})}[e_i]_{\sim},
\end{align}
and $\lim_{i\to \infty}\mu_i=0$. Here and in the following, we write $[e_i]_{\sim}$ if we refer to the $L^2(\nu_{\mc{X}})$ equivalence class and $e_i$ if we refer to the representative in $\mc{H}$.

Further, if there exists $f\in \mch$ such that $[f]_\sim = f^{\ast}$ in $L^2(\nu_\mcx)$, we write by an abuse of notation $f^{\ast}\in \mch$ and identify $f^{\ast}$ with its $\mch$ representative of minimal $\mch$-norm.

Combining the representation of $f^\lambda$ provided by \eqref{eq_min_ridge_general} and \eqref{eq_T_k_L2_classical_1}, we obtain the following classical result for the approximation error.
\begin{theorem}\label{theorem_approximation_error}
    Assume that Assumption \ref{ass_kernel_diagonal} holds.
    Then, for all $\lambda>0$ it holds that
    \begin{align}\label{eq_L2_error_estimate}
        \big\|f^{\ast}-f^{\lambda}\big\|_{L^2(\nu_{\mc{X}})}^2
        =
        \sum_{i\in \mb{N}}\frac{\lambda^2}{(\mu_i+\lambda)^2}\langle f^{\ast},[e_i]_{\sim}\rangle^2_{L^2(\nu_{\mc{X}})},
    \end{align}
    where $\{e_i\}_{i\in \mb{N}}\subset \mch$ is the $\mc{H}$ orthogonal system of eigenvectors with eigenvalues $\{\mu_i\}_{i\in \mb{N}}$ from \eqref{eq_T_k_L2_classical_1}.
    Further, provided we have $f^{\ast}\in \mch$, we find
    \begin{align}\label{eq_H_error_estimate}
        & \|f^{\ast}-f^\lambda\|^2_{\mch}
        =
        \sum_{i\in \mb{N}} \frac{\lambda^2}{\mu_i(\mu_i+\lambda)^2}
        \langle f^{\ast},[e_i]_{\sim}\rangle^2_{L^2(\nu_{\mc{X}})}.
    \end{align}
\end{theorem}
The proof of Theorem \ref{theorem_approximation_error} is standard, see \cite[Theorem 3.5]{BayerBelomestnyButkovskyEtAl2024}. For completeness, we include it in Appendix \ref{sec:approx-error:proofs}.
\begin{remark}
    Observe that since for all $c> 0$ the map $\lambda \mapsto \lambda^2/(c+\lambda)^2$ is monotone, Theorem \ref{theorem_approximation_error} implies that the quantities in \eqref{eq_L2_error_estimate} and \eqref{eq_H_error_estimate} vanish as $\lambda \to 0$.
\end{remark}
In general, the eigenvalues and eigenfunctions of the operator $\mc{T}$ are intractable. Still, when considering  instead of the geometric rough path signature, i.e., Stratonovich, the signature defined by Itô integration, the eigenvalues and eigenfunctions for the case of Brownian motion without time augmentation can be computed exactly. 
As the RKHS $\mch$ is then not dense in $L^2(\nu_\mcx)$ anymore, we get an additional bias term. 
Denoting the $L^2(\nu_\mcx)$ closure of $\mch$ by $\bar{\mch}$, this term measures the distance between  
$f^{\ast}$ and its $L^2(\nu_\mcx)$ projection on $\bar{\mch}$, denoted by $P_{\bar{\mch}} f^{\ast}$.
We obtain an explicit formula for the approximation error, which is the content of the next result. 
\begin{corollary}\label{corollary:approx-error:Brownian}
    Consider the setup where the Stratonovich signature is replaced by the Itô signature.
    If $\nu_{\mc{X}}$ is the $\mb{R}^d$-Wiener measure with time horizon $[0,T]$ and $\nu_{\mb{R}}\in \mc{P}_2(\mb{R})$,
    it holds that
    \begin{align}\label{eq:approx-error:Brownian}
        \big\|f^{\ast}-f^{\lambda}\big\|_{L^2(\nu_{\mc{X}})}^2
        =&
        \big\|f^{\ast}-P_{\bar{\mch}} f^{\ast}\big\|_{L^2(\nu_{\mc{X}})}^2
        \\ \notag &+
        \sum_{m\geq 0}
        \frac{m!}{T^m}
        \frac{\lambda^2}{(T^m/m!+\lambda)^2}
        \sum_{w\in W_m}\langle f^{\ast},\sig^w\rangle^2_{L^2(\nu_{\mc{X}})}.
    \end{align}
\end{corollary}
\begin{remark}
    In the case of Itô signatures for Brownian motion without time augmentation, the entries of the signature are an $L^2$ orthogonal system of $\mc{T}$, due to their $L^2$ orthogonality. Both when adding time augmentation or when switching to Stratonovich signatures, this property is lost. In particular, an explicit computation of \eqref{eq:approx-error:Brownian} seems to be intractable.
\end{remark}
\subsection{Finite-sample analysis under a Bernstein condition}\label{sec:finite-sample:bernstein}
In this section, we present our results for the finite-sample error under a Bernstein condition on the signature kernel, following the approach used by Caponnetto and de Vito \cite{Caponnetto2007OptimalRates}. The proofs can be found in Section \ref{section_sub_finite_sample_error}.
Our goal is to obtain a quantitative bound on the error $f^\lambda_N\approx f^{\ast}$. Given that we have studied already the approximation error $f^\lambda\approx f^{\ast}$ in Section \ref{sec:approx-error}, what remains to do is to analyze the statistical error $f^\lambda_N \approx f^\lambda$.

Introducing for $f\in L^2(\nu_{\mc{X}})$ the expected risk
\begin{align}
    \mc{E}[f]\coloneqq \int_{\mcx \times \R}\,|f(\bx)-y|^2 \,\nu(\mathrm{d}\bx,\mathrm{d}y),
\end{align}
we consider the excess risk error criterion
\begin{align}\label{eq_error_diff}
    \mc{E}_{\textnormal{rel}} [f]
    \coloneqq \mc{E}[f]-\mc{E}[f^{\ast}].
\end{align}
Observe that \eqref{eq_error_diff} is a relative error criterion, in the sense that it measures how well $f_{N}^\lambda$ approximates the true relation between $\bx$ and $y$ as encoded by $\nu$, relative to the $\mc{E}$-optimal function $f^{\ast}$.

We will establish under what scaling of $\lambda$ and $N$ the (random) error criterion $\mc{E}_{\textnormal{rel}} [f^\lambda_N]$ vanishes with high probability as the number of samples $N$ tends to infinity.

For this, we need to impose further requirements on the measure $\nu$.
Importantly, different to a lot of known results in the literature about kernel ridge regression \cite{Caponnetto2007OptimalRates, steinwart_KRR}, we will not directly impose boundedness assumptions on the kernel, but only assume Bernstein-type conditions on the $\nu_{\mcx}$-moment growth of it. In Section \ref{sec:finite-sample:fuk-nagaev}, we will consider more general finite moment assumptions.
Further, we will impose an assumption that requires $\bx$ to be sufficiently informative about $y$. 
\begin{assumption}\label{assumption:finitesample:bernstein}
     The following two properties hold.
    \begin{itemize}
        \item There exists $C_{\nu}>0$ such that for all $m\in \mb{N}$ we have
        \begin{align}\label{eq:ass:kernelmoments}
            \int_{\mc{X}}k(\bx,\bx)^m\,\nu_{\mc{X}}(\mathrm{d}\bx)
            \leq m!\,C_{\nu}^m.
        \end{align}
        \item There exist $\Sigma, M>0$ such that for $\nu_\mcx$-almost every $\bx\in \mcx$ and $m\geq 2$ we have
        \begin{align}\label{eq:ass:noise:bernstein}
            \int_{\mb{R}}|y-f^{\ast} (\bx)|^m\,\nu(\mathrm{d}y|\bx) \leq \frac{1}{2} (m/2)!\,\Sigma^2M^{m-2}.
        \end{align}
    \end{itemize}
\end{assumption}
\begin{remark}
    The condition \eqref{eq:ass:noise:bernstein} is stronger than what is needed for the arguments in \cite{Caponnetto2007OptimalRates}, where an $m!$ instead of an $(m/2)!$ growth bound suffices.
\end{remark}
Under Assumption \ref{assumption:finitesample:bernstein}, similar to \cite[Theorem 4]{Caponnetto2007OptimalRates}, we obtain a bound in probability on $\mc{E}_{\textnormal{rel}}[f_{N}^\lambda]$. The bound we obtain is parametrized by 
\begin{enumerate}
    \item the \textbf{residual},
    \begin{align}
        \label{eq:residual}
        \mc{A}(\lambda) \coloneqq 
        \mc{E}_{\textnormal{rel}} [f^\lambda]
        =
        \mc{E}[f^\lambda]-\mc{E}[f^{\ast}],
    \end{align}
    \item and the \textbf{reconstruction error},
    \begin{align}\label{eq_reconstruction_error}
        \mc{B}(\lambda) \coloneqq \|f^\lambda-f^{\ast}\|_{\mch}^2.
    \end{align}
\end{enumerate}
\begin{remark}
    Note that different to the results in \cite{Caponnetto2007OptimalRates}, 
    our bounds do not depend on
    the effective dimension $\mc{N}(\lambda)= \textnormal{Tr}\left(\frac{\mc{T}}{\mc{T}+ \lambda}\right)$,
    which is a summary statistic for how many eigenvectors of $\mc{T}$ are 
    active at regularization level $\lambda$.
    The issue is that the lack of boundedness of the signature kernel does
    not allow us to isolate the problem of bounding $k$ from bounding other
    quantities in integral equations, which would otherwise lead
    to the effective dimension appearing in the bound.
\end{remark}
The following is our bound on the excess risk in terms of $\mc{A}(\lambda),\mc{B}(\lambda)$.
\begin{theorem}\label{theorem:excess-risk:Bernstein}
    Let $N\in \mb{N},\lambda>0,\eta \in (0,1]$, assume that Assumption \ref{assumption:finitesample:bernstein} holds, and $f^{\ast}\in \mch$. 
    Then, with probability at least $1-\eta$, it holds that
    \begin{align}\notag
        \Erel [f_{N}^\lambda]
        \le & 
        3\mc{A}(\lambda) + \left(C_\nu \log(6/\eta)\frac{1}{\lambda N}\right)\cdot(96 \Sigma^2 + 108 C_\nu \mc{B}(\lambda)) 
        \\ \label{eq:bound:bernstein} & + C_\nu \log^2(6/\eta) \frac{1}{\lambda N^2}\cdot(96M^2 + 96 C_\nu \mc{B}(\lambda)),
    \end{align}
    provided $\lambda, N$ fulfill
    \begin{align}\label{eq_condition_params}
        N\geq
        \left(272 \frac{C_\nu^2}{\lambda^2} \vee 1\right)
        \log\frac{6}{\eta},
    \end{align}
    where $C_\nu>0$ is from Assumption \ref{assumption:finitesample:bernstein}.
\end{theorem}
\begin{remark}\label{remark:explanation}
    Observe that in the above Theorem \ref{theorem:excess-risk:Bernstein}, the function $f_{N}^\lambda$ is a random function, as it is a function of the training data. The $1-\eta$ probability is thus w.r.t. the function $f_{N}^\lambda$, i.e., with probability at least $1-\eta$, the random $N$ i.i.d. training samples are such that $f_{N}^\lambda$ fulfills the above bound \eqref{eq:bound:bernstein}.
\end{remark}
The main tool to prove Theorem \ref{theorem:excess-risk:Bernstein} is a Bernstein-type concentration inequality for random variables in Hilbert spaces \cite[Theorem 3.3.4]{Yurinsky1995}, see also Section \ref{sec:concentration} in this work.

The next Theorem contains a quantitative bound on the excess risk for a fixed $\lambda=\lambda_N$, provided $f^{\ast}\in \mc{H}$. It is our main result for the case $f^{\ast}\in \mc{H}$. It follows as a direct consequence of Theorem \ref{theorem:excess-risk:Bernstein} and explicit computations of the quantities $\mc{A}(\lambda)$ and $\mc{B}(\lambda)$.
\begin{theorem}\label{theorem:rates:Bernstein}
    In the setup of Theorem \ref{theorem:excess-risk:Bernstein}, assume that $f^{\ast}\in \mch$ and let $\eta\in (0,1]$. Then, setting
    \begin{align}
        \lambda_N
        \coloneqq &
        \sqrt{\frac{\log(6/\eta)}{N}}
        \left(6\frac{\sqrt{C_\nu}}{\|f^{\ast}\|_\mch}
        (\Sigma^2 + M^2 \log(6/\eta)/N)^{1/2}
        + 17C_\nu \right)
    \end{align}
    we have with probability at least $1-\eta$ that
    \begin{align}
        \Erel [f_{N}^{\lambda_N}]
        \leq 
        63
        \sqrt{\frac{\log (6/\eta)}{N}} \sqrt{C_\nu}
        \|f^{\ast}\|_\mch
        ((\Sigma^2 + M^2 \log (6/\eta) /N )^{\frac{1}{2}}
        + \sqrt{C_\nu}
        \|f^{\ast}\|_\mch
        ),
    \end{align}
    provided $N\geq \log(6/\eta)$.
\end{theorem}
Observe that in the above we set $\lambda_N=\infty$ if $f^{\ast}=0$, which enforces $f^\lambda=f^\lambda_N=0$. Of course, this is not practical, as there is no need for doing regression if we know a priori that $f^{\ast}=0$.

When not assuming $f^{\ast}\in \mc{H}$, we still obtain the following consistency guarantee as a consequence of an adjusted version of Theorem \ref{theorem:excess-risk:Bernstein}.
\begin{theorem}\label{theorem_consistency}
    Assume that Assumption \ref{assumption:finitesample:bernstein} holds
    and that there is $\lambda_0>0$ such that $\sup_{0<\lambda \leq \lambda_0}\|f^\lambda -f^{\ast}\|_{L^6(\nu_\mcx)}<\infty$. Then, for any sequence $\{\lambda_N\}_{N\in \mb{N}}$ fulfilling 
    \begin{align}\label{eq_condition_no_rates}
        \lim_{N\to \infty} \lambda_N^2 N\wedge 1/\lambda_N=\infty,
    \end{align}
    we have 
    $\lim_{N\to \infty} \Erel[f_{N}^{\lambda_N}]\overset{\mb{P}}{=}0$.
\end{theorem}
Since the Bernstein conditions in Assumption \ref{assumption:finitesample:bernstein} are strictly stronger than the heavy tailed conditions in Assumption \ref{assumption:finitesample:heavytail} stated in the next section, consistency for heavy tailed kernels implies consistency for kernels fulfilling a Bernstein condition with the same scaling in $\lambda, N$. Hence, we omit the proof of consistency for the Bernstein setting, i.e., Theorem \ref{theorem_consistency}, and refer to the proof for the analogous result in the heavy tailed case, i.e., Theorem \ref{theorem:consistency:heavytail}.
\subsection{Finite-sample analysis for heavy tailed signature kernels}\label{sec:finite-sample:fuk-nagaev}
In this section, we discuss convergence guarantees for SKRR when the kernel is heavy tailed, i.e., the Bernstein-type moment growth of the signature kernel from Assumption \ref{assumption:finitesample:bernstein} does not hold. The proofs are contained in Section \ref{sec:heavytailed:proofs}, except for the proof of Lemma \ref{lemma:BM-kernel-moments}, which is in Section \ref{sec:sig-kernel-moment:BM}.

To motivate the interest in this setting, we start by discussing the moment growth of signature kernels with time-augmented Brownian data. The next Lemma provides an upper bound.
\begin{lemma}\label{lemma:BM-kernel-moments}
    Denote by $\bar{B}$ a time augmented Brownian motion on $\mb{R}^{d+1}$ with time horizon $[0,T]$.
    Then, if we consider the signature kernel $k$, there exists a universal constant $C^\ast>0$ such that for all $m\in \mb{N}$ it holds that
    \begin{equation*}
        \mb{E}[k(\bar{B},\bar{B})^m]
        \leq C^\ast \cdot
        2^{m-1}
        \left(1+
        \frac{\sqrt{m}}{(2\pi)^m}
        \exp(4m^2 e (d+1)(T\vee 1)^2+ m)\right).
    \end{equation*}
\end{lemma}
Observe that the above bound scales like $e^{c\cdot m^2}\in \omega(m!e^{c\cdot m})$, so it does not ensure that Assumption \ref{assumption:finitesample:bernstein} holds. In fact, the bound of Lemma \ref{lemma:BM-kernel-moments} is tight, thus Assumption \ref{assumption:finitesample:bernstein} provably does not hold.

To verify this, consider the case $d=1$ and BM without time augmentation and with Stratonovich signature, i.e., the geometric rough path signature. This case bounds below the moments of all other Brownian Stratonovich signature kernels, with and without time augmentation. In that case, we have that
\begin{align}\notag
    \mb{E}\left[k(B,B)^m\right]
    \geq &
    \mb{E}\left[\left(\sum_{n \in \mb{N}} \frac{B_T^{2n}}{(2n)!}\right)^m\right] 
    =
    \mb{E}\left[
    \textnormal{cosh}(B_T)^m
    \right] 
    \geq
    \mb{E}\left[
    \frac{\exp \left( m B_T \right)}{2^m}
    \right]
    .
\end{align} 
Since $\mb{E}\left[\exp\left(mB_T\right)\right]=\exp(m^2 T/2)$, we find
\begin{align}
  \label{eq:lb-moments-strat-bm}
    \mb{E}\left[k(B,B)^m\right]
    \geq 
    \exp(m^2 T/2- m \log 2)
    \gtrsim
    \exp(m^2 T/4).
\end{align}
In particular, we have that $\mb{E}\left[k(B,B)^m\right]$ is of order $e^{c\cdot m^2}$, so Assumption \ref{assumption:finitesample:bernstein} is not fulfilled.

Following the same approach as for Stratonovich signatures, one can prove analogous results for Itô signatures, i.e., also in that case $\mb{E}[k(\bar{B},\bar{B})^m]$ is of order $e^{c\cdot m^2}$.
To do so, one uses that Itô iterated integrals of time augmented Brownian motion are elements of the Wiener chaos, whose moments can be bounded by Nelson's hypercontractivity estimates \cite{Nelson1973FreeMarkoffField, hypercontractivity_remark}. Also the $1d$-approach from the Stratonovich case for the lower bound applies, using the Hermite-polynomial representation for iterated integrals of Brownian motion.

Motivated by Assumption \ref{assumption:finitesample:bernstein} failing to hold in the Brownian case, we now consider instead the following assumption, where we assume the signature kernel to admit only a finite number of moments $m\geq 3$. 
\begin{assumption}\label{assumption:finitesample:heavytail}
     For $m\in \N,m\geq 3$, the following two properties hold.
    \begin{itemize}
        \item There exist constants $(C_{m'})_{1\leq m'\leq m}\subset (0,\infty)$, such that for $1\leq m' \leq m$ we have
        \begin{align}\label{eq:ass:kernelmoments:heavytailed}
            \int_{\mc{X}}k^{m'}(\bx,\bx)\,\nu_{\mc{X}}(\mathrm{d}\bx)
            \leq C_{m'}.
        \end{align}
        \item 
        There exist constants $(D_{m'})_{1\leq m'\leq m}\subset (0,\infty)$, such that for $1\leq m' \leq m$ and $\nu_\mcx$-almost every 
        $\bx\in \mcx$ we have
        \begin{align}\label{eq:assumption:heavytail:y-f}
            \int_{\mb{R}}|y-f^{\ast} (\bx)|^{m'}\,\nu(\mathrm{d}y|\bx) \leq D_{m'}.
        \end{align}
    \end{itemize}
\end{assumption}
Observe that by Lemma \ref{lemma:BM-kernel-moments}, the first part of Assumption \ref{assumption:finitesample:heavytail} is fulfilled if $\nu_\mcx$ is the law of Brownian motion with time augmentation.
Further, even though Assumption \ref{assumption:finitesample:heavytail} does not allow anymore for applications of Bernstein's inequality, one can still replace Bernstein by the weaker Fuk-Nagaev inequality, which we introduce in Section \ref{sec:concentration}.

The results for the heavy tailed kernel get established in the same way as for the case of Assumption \ref{assumption:finitesample:bernstein}, i.e., using the Fuk-Nagaev inequality, we prove results analogous to Theorems
\ref{theorem:excess-risk:Bernstein}, \ref{theorem:rates:Bernstein}, and \ref{theorem_consistency}.

The following result is analogous to Theorem \ref{theorem:excess-risk:Bernstein}.
\begin{theorem}\label{theorem:errorbound:heavytail}
    Let $N\in \mb{N},\lambda>0,\eta \in (0,1]$, assume that Assumption \ref{assumption:finitesample:heavytail} holds, and $f^{\ast}\in \mch$. 
    Then, with probability at least $1-\eta$, it holds that 
    \begin{align}\notag
            \Erel[f_{N}^\lambda]
            \leq & C'_m \bigg( \mc{A} (\lambda) +\left( \frac{1}{\lambda}
            +\frac{{\mc{B}(\lambda)}}{\lambda}\right)
            \max\left\{
            \left(\frac{1 }{\eta N^{m-1}}\right)^{1/m},
            \sqrt{ \frac{\log (6c_1/\eta)}{N}}
            \right\}^2 
            \bigg),
    \end{align}
    provided $\lambda, N$ fulfill
    \begin{align}\label{condition:errorbound:heavytail}
        \lambda
            \min \left\{ 
            ((\eta/3) N^{m-1})^{1/m},
            \sqrt{N} \log (6c_1/\eta)^{-1/2}
            \right\} \geq C''_m,
    \end{align}
    where the constants $C'_m,C''_m,c_1>0$ depend only on $m$ and the constants from Assumption \ref{assumption:finitesample:heavytail}.
\end{theorem}
As a consequence of the above, we obtain rates in the case where $f^{\ast}\in \mch$ and by a modification of the above result for the case $f^{\ast}\notin \mch$ also consistency in that case, which is the content of the following two Theorems.
\begin{theorem}\label{theorem:rates:heavytail}
    In the setup of Theorem \ref{theorem:errorbound:heavytail}, let $\eta \in (0,1]$. Then, setting 
    \begin{align}
        \lambda_N \coloneqq 
        \frac{\sqrt{\log (6c_1/\eta)}}{\sqrt{N}}
        \left(
        C''_m
        \vee
        \left( \frac{1}{\|f^{\ast}\|_{\mch}^2}
        +1\right)^{1/2}
        \right),
    \end{align}
    with probability at least $1-\eta$, we have
    \begin{align}\notag
        \Erel[f_{N}^{\lambda_N}]
        \leq & C'''_m
        \frac{\sqrt{\log (6c_1/\eta)}}{\sqrt{N}}
        \bigg( 
        \|f^{\ast}\|_{\mch}^2 
        \left(
        C''_m
        \vee
        \left( \frac{1}{\|f^{\ast}\|_{\mch}^2}
        +1\right)^{1/2}
        \right)
        \\ \notag
        & +
        \frac{1}{C_m''}
        \left( 1
        +\|f^{\ast}\|_{\mch}^2\right)
        \bigg),
    \end{align}
    provided $N$ fulfills additionally $N \geq (\eta/3)^{-2/(m-2)} 
    \log (6c_1/\eta)^{-m/(m-2)}$. Here, $C'''_m>0$ depends only on $m$ and the constants from Assumption \ref{assumption:finitesample:heavytail}.
\end{theorem}
\begin{remark}\label{remark:rates:comparison}
    Observe that both under the Bernstein and the Fuk-Nagaev assumption, the 
    bounds on the excess risk obtained in Theorems \ref{theorem:rates:Bernstein} and \ref{theorem:rates:heavytail} scale like $1/\sqrt{N}$. Besides the different constants involved in the bounds, a key difference in the results is that in the Fuk-Nagaev setting, one needs to impose a stronger lower bound on $N$ in terms of $\eta$.
\end{remark}
The next Theorem derives consistency under the same scaling as in Theorem \ref{theorem_consistency}.
\begin{theorem}\label{theorem:consistency:heavytail}
    Assume that Assumption \ref{assumption:finitesample:heavytail} holds
    and that there is $\lambda_0>0$ such that $\sup_{0<\lambda \leq \lambda_0}\|f^\lambda -f^{\ast}\|_{L^6(\nu_\mcx)}<\infty$. Then, for any sequence $\{\lambda_N\}_{N\in \mb{N}}$ fulfilling 
    \begin{align}\label{eq_condition_no_rates:heavytail}
        \lim_{N\to \infty} \lambda_N^2 N\wedge 1/\lambda_N = \infty,
    \end{align}
    we have 
    $\lim_{N\to \infty} \Erel[f_{N}^{\lambda_N}]\overset{\mb{P}}{=}0$.
\end{theorem}
\subsection{Hilbert space concentration inequalities}\label{sec:concentration}
The main tools to establish the previously stated results are concentration inequalities for random variables taking values in a Hilbert space, in our case the RKHS $\mch$. One of the most classical concentration inequalities is due to Bernstein, for which one assumes the random variable to admit sub-exponential moments. The version of Bernstein in a Hilbert space setting stated next is taken from \cite{Yurinsky1995}.
\begin{theorem}[{\cite[Theorem 3.3.4]{Yurinsky1995}}]\label{theorem:Bernstein}
    Consider a random variable $\zeta$ taking values in a real separable Hilbert space $\mc{K}$. Assume that there are two positive constants $L,\sigma>0$ such that
    \begin{align}
        \mb{E}\big[\|\zeta-\mb{E}[\zeta]\big\|^m_{\mc{K}}\big]
        \leq \frac{1}{2}m! \sigma^2 L^{m-2}\,\,\,\,\,\,\,\,\forall m\geq 2.
    \end{align}
    Then, denoting for $N\in \mb{N}$ by $\zeta_1,\dots,\zeta_N$ i.i.d. copies of $\zeta$, for all $\eta\in (0,1)$ it holds that
    \begin{align}\notag
        \mb{P}\left(
        \left\|\frac{1}{N}\sum_{1\leq i \leq N} \zeta_i 
        -\mb{E}[\zeta]\right\|_{\mc{K}}
        \leq 
        \frac{2 L}{N} \log \frac{2}{\eta}
        +\sigma
        \sqrt{\frac{2}{N}
        \log \frac{2}{\eta}}
        \right)
        \geq 1 - \eta.
    \end{align}
\end{theorem}
Next, we present the Fuk-Nagaev inequality. The version stated is a modification of \cite[Theorem 3.5.1]{Yurinsky1995}, that is due to \cite{mollenhauer2025regularized}. For a discussion of the involved constants, see \cite{Rio2017}.
\begin{theorem}[{\cite[Corollary 5.3]{mollenhauer2025regularized}}]\label{theorem:fuknagaev}
    Consider a random variable $\zeta$ taking values in a real separable Hilbert space $\mc{K}$. Let $m\in \N,\,m\geq 3$ and assume that there are two positive constants $L_m,\sigma>0$ such that
    \begin{align}
        \mb{E}\big[\|\zeta-\mb{E}[\zeta]\big\|^2_{\mc{K}}\big]
        \leq \sigma^2\quad \textnormal{and}\quad
        \mb{E}\big[\|\zeta-\mb{E}[\zeta]\big\|^m_{\mc{K}}\big]
        \leq L_m.
    \end{align}
    Then, denoting for $N\in \mb{N}$ by $\zeta_1,\dots,\zeta_N$ i.i.d. copies of $\zeta$, for all $\eta\in (0,1)$ it holds that
    \begin{align}\notag
        \mb{P}\left(
        \left\|\frac{1}{N}\sum_{1\leq i \leq N} \zeta_i
        -\mb{E}[\zeta]\right\|_{\mc{K}}
        \leq 
        \max\left\{
        \left(\frac{2c_1 L_m }{\eta N^{m-1}}\right)^{1/m},
        \sigma \sqrt{ \frac{\log (2c_1/\eta)}{c_2 N}}
        \right\}
        \right)
        \geq 1 - \eta,
    \end{align}
    where the constants $c_1,c_2>0$ depend only on $m$.
\end{theorem}

\section{Robust signature kernel ridge regression}
In this section, we discuss \textit{robust} signature kernel ridge regression, where the signature is replaced by a properly normalized signature, referred to as the robust signature, introduced by Chevyrev and Oberhauser \cite{chevyrev_oberhauser_robust}. The motivation for considering this map 
is that it has better statistical properties than the standard signature. 
It is constructed by injectively mapping signatures to the centered ball of radius $\sqrt{\kappa}>1$ in $T((\R^{d+1}))$, in such a way that the obtained robust signature $\widehat{\sig}(\bx)$ is the signature of a rescaling of $\bx$. The robust signature kernel is then bounded by $\kappa$.
Also, analogously to the standard case, the RKHS $\widehat{\mch}$ of the robust signature kernel is a space of linear functionals of the robust signature. 
For more details, we refer to Appendix \ref{appendix:robust}.

The approach for deriving theoretical guarantees for robust SKRR is the same as for classical SKRR. The following assumption is analogous to Assumption \ref{ass_kernel_diagonal} and assumed to hold in the following.
\begin{assumption}\label{ass_kernel_diagonal:robust}
    The measure $\nu\in \mc{P}(\mcx\times \R)$ is such that $\nu_{\R }\in \mc{P}_2(\R)$.
    Further, the RKHS $\widehat{\mch}$ is dense in $L^2(\nu_\mcx)$.
\end{assumption}

The analysis of the approximation error proceeds in the same way as in the standard case and one obtains a result analogous to Theorem \ref{theorem_approximation_error}. Note though that computing eigenvalues and eigenfunctions is more challenging, due to the normalization of the robust signature.

To derive a finite-sample guarantee analogous to Theorem \ref{theorem:rates:Bernstein}, since the robust signature kernel is bounded, we can apply the results of Caponnetto and de Vito \cite{Caponnetto2007OptimalRates}. We work under the following assumption.
\begin{assumption}\label{assumption_nu_finite_sample:robust}
    There exist $\Sigma, M>0$ such that for $\nu_\mcx$-almost every $\bx\in \mcx$ and $m\geq 2$ we have
    \begin{align}
        \int_{\mb{R}}|y-f^{\ast} (\bx)|^m\,\nu(\mathrm{d}y|\bx) \leq \frac{1}{2} m!\,\Sigma^2M^{m-2}.
    \end{align}
\end{assumption}
\begin{remark}
    Observe that due to the boundedness of the robust signature kernel, we do not need to impose any moment growth assumptions. Further, this boundedness also allows us to relax the conditional moment growth assumption on $y-f^{\ast}(\bx)$, i.e., we replaced $(m/2)!$ by $m!$.
    The assumption considered in \cite{Caponnetto2007OptimalRates} is stronger, though for the results of \cite{Caponnetto2007OptimalRates} relevant to us, Assumption \ref{assumption_nu_finite_sample:robust} is sufficient.
\end{remark}
The following result is our finite-sample guarantee.
We indicate the analogues of the objects from the standard signature 
case for the robust case by a hat accent, i.e., the analogues of 
$f^\lambda,f_{N}^\lambda$ are denoted by $\hat{f}^\lambda,\hat{f}_{N}^\lambda$.
\begin{theorem}\label{theorem:rates:robust}
    Assume that Assumption \ref{assumption_nu_finite_sample:robust} holds and that $f^{\ast}\in \mchrob$. Then, for any $\eta \in (0,1]$, setting
    \begin{align}\label{eq:lambdaN:robust}
        \lambda_N \coloneqq
        \frac{1}{\sqrt{N}} 
        \left(
        (\frac{M^2}{308}+\Sigma^2)^{1/2}
        \frac{\sqrt{\Tr (\widehat{\mc{T}})}}{\|f^{\ast}\|_{\mchrob}}
        \vee \sqrt{2 C_{\eta} \kappa \Tr (\widehat{\mc{T}})}
        \right),
    \end{align}
    where $C_\eta = 96 \log^2(6/\eta)$ and $\textnormal{Tr}(\widehat{\mc{T}})$ denotes the operator trace of $\widehat{\mc{T}}$,
    we have with probability at least $1-\eta$ that
    \begin{align}\notag
        \mc{E}_{\textnormal{rel}}[\hat{f}^{\lambda_N}_N]
        \leq  
        \sqrt{\frac{\Tr (\widehat{\mc{T}})}{N}} 
        2
        C_\eta
        \|f^{\ast}\|_{\mchrob}
        \left(
        2
        \left(\frac{M^2}{102}+\Sigma^2\right)
        +
        2 C_\eta \kappa
        \|f^{\ast}\|_{\mchrob}^2
        \right)^{1/2},
    \end{align}
    provided $N$ fulfills
    \begin{align}\label{eq:conditionN:robust}
        N\geq 
        \frac{1}{\|\widehat{\mc{T}}\|_{\mc{L}(\mchrob)}^2}
        \left(
        (\frac{M^2}{102}+\Sigma^2)
        \frac{\Tr (\widehat{\mc{T}})}{\|f^{\ast}\|_{\mchrob}^2}
        \vee
        C_\eta 2 \kappa \Tr (\widehat{\mc{T}}) \right).
    \end{align}
    Here, $\kappa>1$ denotes the bound on the robust signature kernel $\hat{k}$.
\end{theorem}
A proof of Theorem \ref{theorem:rates:robust} can be found in the Appendix.
\begin{remark}\label{remark:comarison-standard-robust}
    Note that also in Theorem \ref{theorem:rates:robust}, the bound on the excess risk scales as $1/\sqrt{N}$. Observe further that, due to the boundedness of the robust signature kernel, the bounds on its moments can be chosen distribution independent, which is in sharp contrast to the results of Theorems \ref{theorem:rates:Bernstein} and \ref{theorem:rates:heavytail}, where explicit moment assumptions are needed.
\end{remark}

In the case where we do not assume $f^{\ast}\in \mchrob$, we get the same consistency result as in the standard case, though under a weaker moment assumption on $\hat{f}^\lambda-f^{\ast}$. 
\begin{theorem}\label{theorem_consistency_robust}
    Assume that Assumption \ref{assumption_nu_finite_sample:robust} holds
    and that there is $\lambda_0>0$ such that $\sup_{0<\lambda \leq \lambda_0}\|\hat{f}^\lambda -f^{\ast}\|_{L^3(\nu_\mcx)}<\infty$.
    Then, for any sequence $\{\lambda_N\}_{N\in \mb{N}}$ fulfilling 
    $$\lim_{N\to \infty} \lambda_N^2 N\wedge 1/\lambda_N=\infty,$$ we have 
    \begin{align}
        \lim_{N\to \infty} 
        \mc{E}_{\textnormal{rel}}[\hat{f}^{\lambda_N}_N]
        \overset{\mb{P}}{=}0.
    \end{align}
\end{theorem}
To prove Theorem \ref{theorem_consistency_robust}, one can apply the whole approach for the standard signature kernel directly to the robust case. Indeed, since the robust signature kernel is bounded and since we require Assumption \ref{assumption_nu_finite_sample:robust}, we are in the framework of the heavy-tailed proof, up to the weaker moment assumption on $\hat{f}^\lambda -f^{\ast}$. We can weaken this assumption from $L^6(\nu_\mcx)$ to $L^3(\nu_\mcx)$ due to the boundedness of the robust signature kernel. Hence, we consider the result as proven and refer to the proof of Theorem \ref{theorem:consistency:heavytail}.

\section{Numerical Experiments}
\label{sec:numer-exper}

In this section, we present results of numerical experiments on two synthetic datasets and one real world dataset with various types of signature kernels. 

For the synthetic experiments, we consider 
training samples of sizes $N_{\textnormal{train}}\in \{2^4,2^5,\dots, 2^9\}$ with fixed test size $N_{\textnormal{test}}=200$. Given training and test samples 
$\{(\bx_i,y_i)\}_{1\leq i \leq N_{\textnormal{train}}}$, $\{(\bx'_i,y'_i)\}_{1\leq i \leq N_{\textnormal{test}}}$, we compute 
\begin{equation*}
    f^\lambda_{N_{\textnormal{train}}} = \argmin_{f\in \mch} 
    \left\{
    \frac{1}{N_{\textnormal{train}}} \sum_{1 \leq i \leq N_{\textnormal{train}}} \left(f(\bx_i)-y_i\right)^2 + \lambda \|f\|^2_{\mch}
    \right\}
    ,
\end{equation*}
and then evaluate $f^\lambda_{N_{\textnormal{train}}}$ by computing
\begin{equation*}
    \frac{1}{N_{\textnormal{test}}}
    \sum_{1 \leq i \leq N_{\textnormal{test}}} (f^\lambda_{N_\textnormal{train}}(\bx'_i)-f^{\ast}(\bx'_i))^2,
\end{equation*}
i.e., we compute excess risk. The estimator is justified
by the identity $\mc{E}_{\textnormal{rel}}[f]= \int_{\mcx} (f(\bx)-f^{\ast}(\bx))^2 \,\nu_{\mcx}(\dd \bx)$.
We run each of the experiments $N_{\textnormal{MC}}=50$ times and then compute an estimator of the losses as the average (for any fixed $\lambda$ and $N_{\textnormal{train}}$) across the $N_{\textnormal{MC}}=50$ repetitions. We always choose a $\lambda$ grid adapted to the given setup. Note that in our synthetic examples, the conditional expectations $f^{\ast}$ are available in closed form. 

In our real
world time series example, we compute the root mean squared error (RMSE)
\begin{equation*}
    \sqrt{\frac{1}{N_{\textnormal{test}}}
    \sum_{1 \leq i \leq N_{\textnormal{test}}} (f^\lambda_{N_\textnormal{train}}(\bx'_i)-y_i')^2},
\end{equation*}
as the conditional 
expectation is not available. We also use a train-test split, with split size depending on the time series considered. 

In both the synthetic as well as the real-world setup, all paths are augmented by a time component.

We compute the signature kernel either via truncated    signatures or by the PDE approach from \cite{goursat_PDE}. In the latter case, we use the implementation provided by the \texttt{sigkernel} library. For (discrete) $\R^d$ valued paths $\mathbf{x},\mathbf{x}'$ of lengths $L,L'$, the computational complexity of computing $k(\mathbf{x},\mathbf{x}')$ without parallelization using \texttt{sigkernel} and a fixed grid is of order $\mc{O}(d L L')$.
We perform all computations of (truncated) signatures using the $\texttt{iisignature}$ library \cite{reizenstein2018iisignaturelibraryefficientcalculation}. It provides efficient algorithms to compute truncated signatures with computational complexity $\mathcal{O}(Ld^m)$, where $m$ denotes the truncation level. 

For robust signatures, we choose the version
of \cite[Example 4]{chevyrev_oberhauser_robust}
and compute them with the approach given in the same work, using the code from \cite{bayer2026local}.

To compute the inverse of $K+N\lambda \textnormal{Id}$, we compute a spectral decomposition of $K$. While this is more computationally intensive for a given $\lambda$ than just inverting $K+N\lambda \textnormal{Id}$, we have to compute the spectral decomposition only once and can then apply it across the whole
grid of $\lambda$ values.

\subsection{One-dimensional SDE ({\cite[section 4.1]{bayer2026local}})}
\label{sec:one-dimensional-sde}

Let $T>0,z_0\in \R^d$, and $B=(B_t)_{0\leq t \leq T}$ be a Brownian motion. Then, the SDE
\begin{align}\label{eq:onedimSDE}
    \begin{cases}
        \mathrm{d}Z_t &= -Z_t^p\,\mathrm{d}t + Z_t\cos(Z_t)\,\mathrm{d}B_t, 
        \\ Z_0 &= z_0
    \end{cases}
\end{align}
admits a unique strong solution on $[0,T]$, see for example \cite[Chapter 3, Theorem 7]{Protter2005}. 

We aim at learning the terminal value of the solution $Z$ given the driving noise $B$ using robust signature kernel ridge regression.

This problem was first considered in \cite[section 4.1]{bayer2026local}, to study a signature-based local regression method. 
Different to \cite{bayer2026local}, to model a noisy setup, we multiply the terminal value $Z_T$ with $1+0.3 \times \chi^2$, where the $\chi^2$ distributed random variable is taken independently of the Brownian motion, that is we observe $(B,(1+0.3 \times \chi^2)Z_T)$ in training.
We use the same parameters as in \cite{bayer2026local}, i.e., we choose $d=1$, $p=5$, and $T=1$. Further, we take $N_{\textnormal{disc}}=400$ steps in the Euler-Maruyama discretization.

\begin{figure}
    \centering
    \begin{minipage}{0.48\textwidth}
        \centering
        \includegraphics[width=\textwidth]{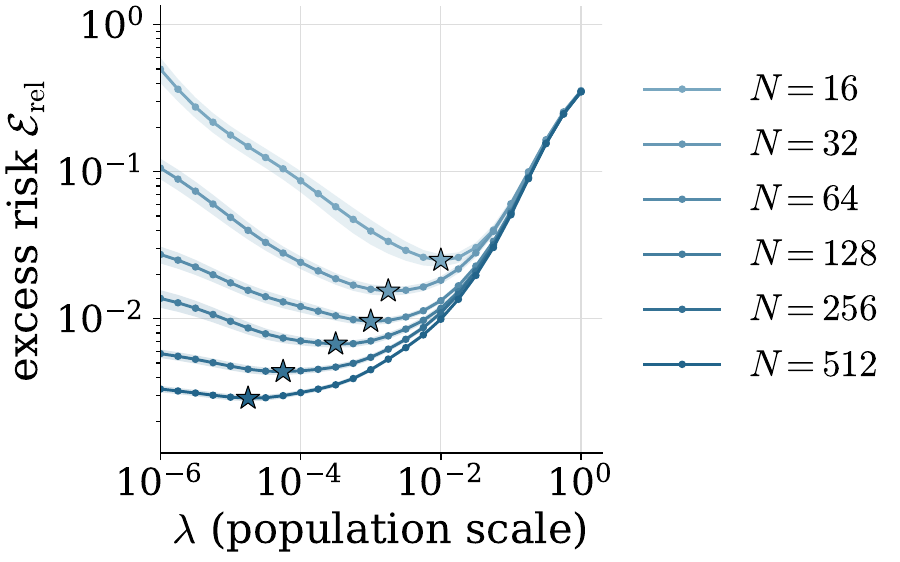}
    \end{minipage}
    \hfill
    \begin{minipage}{0.48\textwidth}
        \centering
        \includegraphics[width=\textwidth]{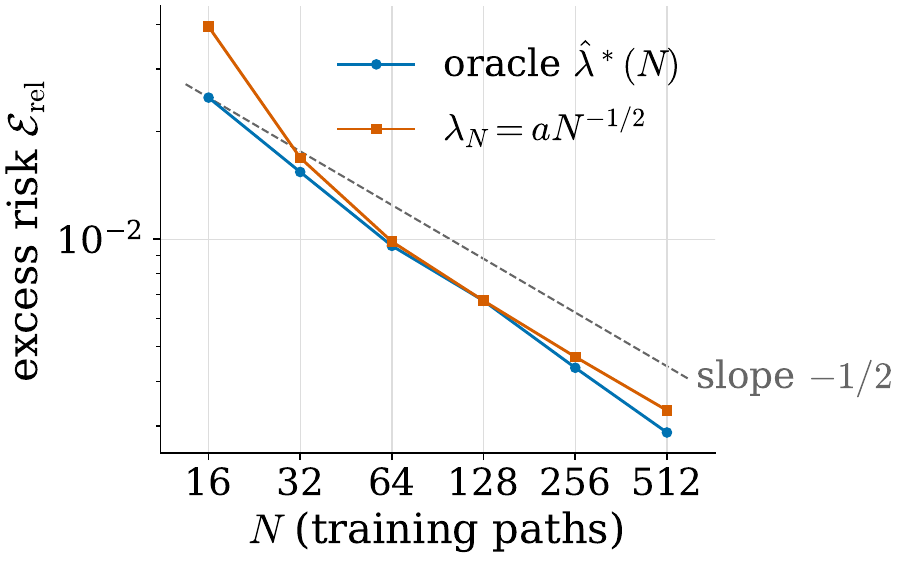}
    \end{minipage}
    \caption{Robust SKRR for $(B_t)_{0\leq t \leq T}\mapsto Z_T \cdot (1+0.3 \chi^2)$, where $Z$ solves \eqref{eq:onedimSDE}. Left: Mean excess risk plotted against $\lambda$ for different training sample sizes $N$. The displayed band corresponds to $\pm$ its standard error at $\lambda,N$. The stars highlight the optimal $\lambda$. Right: Plot of excess risk for optimal $\lambda$ (oracle) and for $\lambda$ approximately fitted to a curve corresponding to the scaling of Theorem \ref{theorem:rates:robust} on log-scales.}
    \label{fig:skrr}
  \end{figure}

Observe that in the left part of Figure \ref{fig:skrr} we see the decay in $\lambda$ as $N$ increases, whereas in the right we see that the $\lambda$ suggested by the scaling of Theorem \ref{theorem:rates:robust} almost achieves the same excess risk as the optimal $\lambda$. 

\subsection{Heston signature model}
\label{sec:hest-sign-model}

\begin{figure}
    \centering
    \begin{minipage}{0.48\textwidth}
        \centering
        \includegraphics[width=\textwidth]{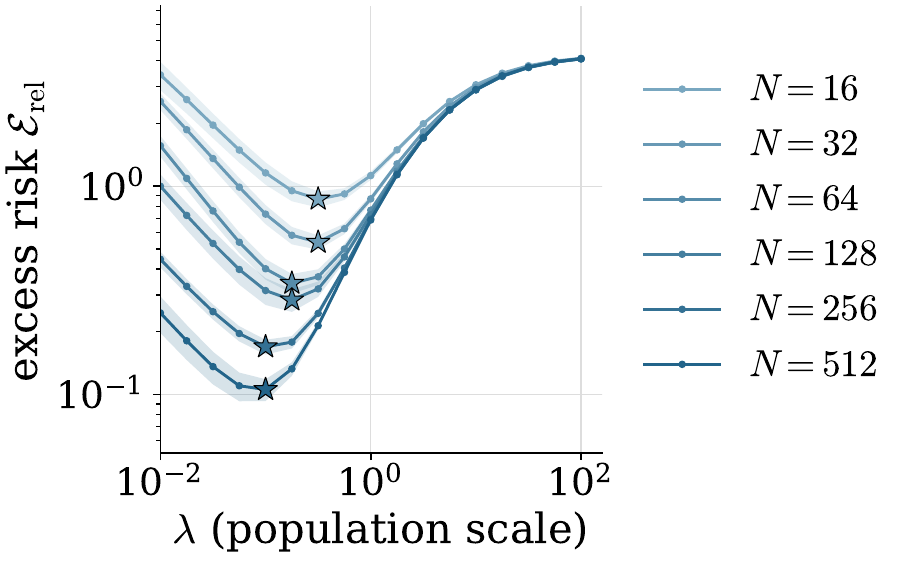}
    \end{minipage}
    \hfill
    \begin{minipage}{0.48\textwidth}
        \centering
        \includegraphics[width=\textwidth]{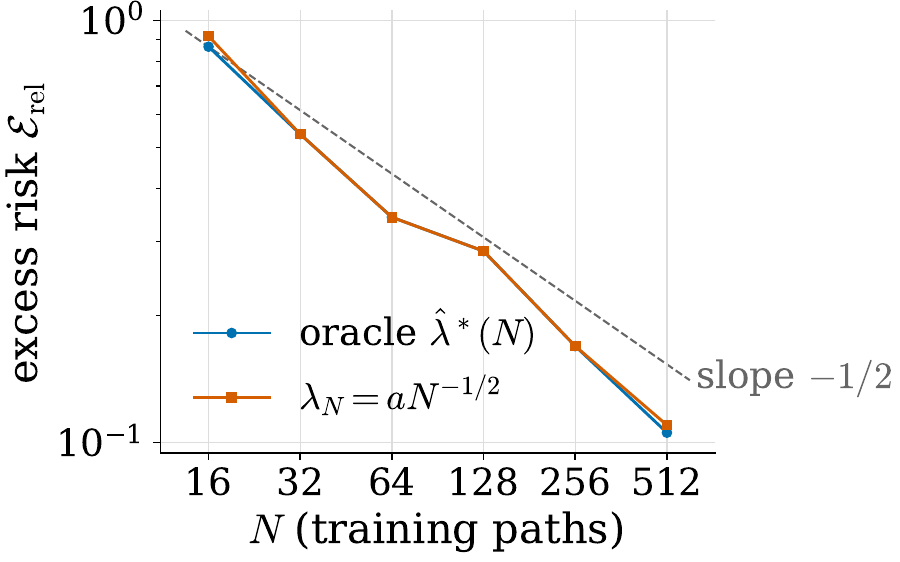}
    \end{minipage}
    \caption{SKRR for a signature approximation of the Heston model. Left: Excess risk plotted against $\lambda$ for different training sample sizes $N$. The displayed band corresponds to $\pm$ its standard error at $\lambda,N$. The stars highlight the optimal $\lambda$. Right: Plot of excess risk for optimal $\lambda$ (oracle) and for $\lambda$ approximately fitted to a curve corresponding to the scaling of Theorem \ref{theorem:rates:heavytail} on log-scales.}
    \label{fig:skrr:Heston}
  \end{figure}

The Heston model is a stochastic volatility model whose dynamics for
stock price and variance are governed by the differential equations
\begin{align}
    & \dd S_t = S_t \mu\,\dd t+ S_t \sqrt{V_t}\,\dd B_t,
    \\ \notag 
    & \dd V_t = \beta (\theta-V_t)\,\dd t + \alpha \sqrt{V_t} \,\dd W_t,
\end{align}
where the Brownian motions $B,W$ driving stock price and variance have 
correlation $\dd \langle B,W\rangle_t = \rho \,\dd t$ with $\rho\in [-1,1]$ and we are given 
initial conditions $S_0,V_0\in \R_+$. We choose $T=1$ with $N_{\textnormal{disc}}=90$ discretization points. 
As done in \cite{cuchiero2022signaturebasedmodelstheorycalibration}, using Lasso we fit a depth three signature model to the Heston model with parameters
$\mu = 0.001,\alpha = 0.25,\beta =0.5,\theta=0.15,\rho=-0.5,S_0=1,V_0 = 0.08$. This model learns, given $X\coloneqq (t,B,W)$ the price $S$, i.e.,
\begin{align}
    \langle \ell , \sig ( X_{[0,t]})\rangle \approx S_t.
\end{align}
Using the learned generative model for Heston prices, we generate data, where we multiply the target with a $\chi^2$ distribution plus a constant to model noise.
We then solve the inverse problem of learning the linear functional of the signature of $X$ using SKRR, given i.i.d. data of $(X,(1+ \chi^2) S_1)$, where $\chi^2$ is independent of the Brownian motions. The signature kernel is computed using truncation with depth $6$. 
Importantly, for this model we have $f^{\ast} \in \mathcal{H}$.
Our results are illustrated in Figure \ref{fig:skrr:Heston}. Observe that its right hand side illustrates decay of the excess risk that is in accordance with
the bound of Theorem \ref{theorem:rates:heavytail}.
\begin{remark}
    Note that the above framework is a priori not covered by Theorem \ref{theorem:rates:heavytail} as stated, since the pointwise condition \eqref{eq:assumption:heavytail:y-f} in Assumption \ref{assumption:finitesample:heavytail} breaks. Nonetheless, since the proof of Theorem \ref{theorem:rates:heavytail} requires only an integrated version of the form
    \begin{equation*}
        \int_{\mcx}
            k(\bx,\bx)^{m/2}
            \left(
            \int_{\R}
            |y - f^{\ast}(\bx)|^{m}
            \, \nu(\mathrm{d}y|\bx)
            \right)
            \, \nu_{\mc{X}}(\mathrm{d}\bx)
            \leq D_m,
    \end{equation*}
    where $D_m>0$, which is fulfilled in our setup, the result of Theorem \ref{theorem:rates:heavytail} holds also in our setup.
\end{remark}

\subsection{TSER repository \cite{TSERTan2020}}
\label{sec:tser}

Time Series Extrinsic Regression \cite{TSERTan2020} is a database of 
regression problems in which one aims at learning the relationship between a time series and a continuous scalar variable.
We interpret the time series as continuous paths by linearly interpolating between the data points.
Since the signature cannot see the initial condition of the time series, we augment every time series with the initial value zero, i.e., every time series starts at zero.

In the preprocessing, we first normalize paths to have zero mean and unit variance on each channel. Then, we normalize paths by a constant to have maximum value equal to one, which is necessary for the Sig-PDE method to be efficiently applicable, see \cite{goursat_PDE}. 

In the numerical experiments, we compute the signature kernel with different truncation levels and via the Sig-PDE approach. For the ApplEnergy dataset we only consider the truncation level $m=2$, as the higher levels are computationally intractable due to the high dimension $d=24$ of the dataset and the exponential scaling of the signature dimension in the truncation level. For the same reason, we do not consider $m=6$ for the BeijingPM25 dataset. We also consider the robust signature kernel with the same choices of truncation levels as in the standard case. For the five datasets with sequence length $\geq 200$, we did not apply Sig-PDE for computational reasons.

The results are illustrated in Tables \ref{table:rmse-lambda} and \ref{table:RMSE-kernel} below, showcasing the benefit of regularization and the different performances of different computational approaches to the standard and robust signature kernel. Observe that in particular on the high dimensional ApplEnergy dataset, where even moderate truncation levels are infeasible, Sig-PDE performs well.
\begin{table}[ht]
\centering
\begin{adjustbox}{max width=\textwidth}
\begin{tabular}{lrrrrrrrrr}
\toprule
Dataset / $\lambda$ & $e^{-11}$ & $e^{-10}$ & $e^{-9}$ & $e^{-8}$ & $e^{-7}$ & $e^{-6}$ & $e^{-5}$ & $e^{-4}$ & $e^{-3}$ \\
\midrule
ApplEnergy & 3.62 & 3.62 & 3.62 & 3.62 & 3.62 & 3.62 & 3.59 & \textbf{3.45} & 3.48 \\
BeijingPM25 & 99.9 & 66.1 & 59.4 & \textbf{55.9} & 59.1 & 64.5 & 67.6 & 70.2 & 86.4 \\
BinanceSent & 0.494 & 0.342 & \textbf{0.318} & 0.322 & 0.328 & 0.335 & 0.34 & 0.34 & 0.348 \\
BitcoinSent & 0.969 & 0.639 & 0.298 & 0.248 & 0.241 & 0.235 & 0.234 & \textbf{0.233} & 0.233 \\
CardanoSent & 1.49 & 0.46 & 0.361 & 0.326 & 0.327 & 0.316 & 0.298 & \textbf{0.297} & 0.306 \\
Covid19And & 0.0389 & 0.0242 & \textbf{0.0229} & 0.0239 & 0.0257 & 0.0264 & 0.0278 & 0.0286 & 0.0295 \\
Covid3M & 0.0438 & \textbf{0.0427} & 0.0427 & 0.0427 & 0.0427 & 0.0427 & 0.0429 & 0.0439 & 0.0496 \\
OilGasPrices & 18.5 & 4.72 & 1.96 & 1.38 & 1.36 & \textbf{1.34} & 1.48 & 1.54 & 1.6 \\
EthereumSent & 0.278 & \textbf{0.224} & 0.224 & 0.224 & 0.224 & 0.225 & 0.226 & 0.229 & 0.26 \\
Flood1 & \textbf{0.0164} & 0.0179 & 0.0183 & 0.0188 & 0.0188 & 0.019 & 0.0351 & 0.188 & 0.388 \\
Flood2 & \textbf{0.0176} & 0.0184 & 0.0185 & 0.0185 & 0.0186 & 0.0191 & 0.0383 & 0.146 & 0.22 \\
Flood3 & \textbf{0.00946} & 0.0166 & 0.0208 & 0.0219 & 0.0226 & 0.0229 & 0.025 & 0.074 & 0.273 \\
GasSensorAc & 0.379 & 0.254 & \textbf{0.25} & 0.267 & 0.272 & 0.271 & 0.272 & 0.276 & 0.282 \\
GasSensorEth & 0.7 & 0.437 & 0.379 & 0.388 & 0.397 & 0.311 & \textbf{0.29} & 0.297 & 0.3 \\
MetroTraffic & 537 & 536 & \textbf{535} & 535 & 535 & 540 & 785 & 2.27e+03 & 3.18e+03 \\
SierraSnow & 73.1 & 36.1 & \textbf{24.8} & 27.4 & 29.8 & 31.9 & 44.6 & 60 & 66.8 \\
\bottomrule
\end{tabular}
\end{adjustbox}
\caption{RMSE for SKRR applied to the TSER dataset across $\lambda$ with the robust signature kernel at truncation level $m=4$, except for the ApplEnergy dataset ($m=2$). Minimal error is highlighted by bold font in each row.}
\label{table:rmse-lambda}
\end{table}

\newpage \begin{table}[ht]
\centering
\label{tab:performance}
\begin{adjustbox}{max width=\textwidth}
\begin{tabular}{lrrrrrrrrr}
\toprule
Dataset & dim & length & SigPDE & RSig(M=2) & RSig(M=4) & RSig(M=6) & Sig(M=2) & Sig(M=4) & Sig(M=6) \\
\midrule
ApplEnergy & 24 & 144 & \textbf{2.87} & 3.45 & -- & -- & 3.35 & -- & -- \\
BeijingPM25 & 6 & 24 & 59.4 & 56.6 & \textbf{55.9} & -- & 56.6 & 55.9 & -- \\
BinanceSent & 2 & 24 & \textbf{0.316} & 0.329 & 0.318 & 0.318 & 0.329 & 0.317 & 0.317 \\
BitcoinSent & 2 & 24 & \textbf{0.231} & 0.233 & 0.233 & 0.233 & 0.233 & 0.233 & 0.233 \\
CardanoSent & 2 & 24 & \textbf{0.296} & 0.297 & 0.297 & 0.297 & 0.298 & 0.298 & 0.298 \\
Covid19And & 1 & 91 & 0.0273 & 0.0273 & 0.0229 & \textbf{0.0215} & 0.0278 & 0.0233 & 0.0233 \\
Covid3M & 1 & 84 & 0.0428 & 0.0427 & 0.0427 & \textbf{0.0426} & 0.0427 & 0.0427 & 0.0427 \\
OilGasPrices & 2 & 30 & \textbf{1.28} & 1.37 & 1.34 & 1.34 & 1.55 & 1.34 & 1.34 \\
EthereumSent & 2 & 24 & \textbf{0.221} & 0.223 & 0.224 & 0.224 & 0.223 & 0.224 & 0.224 \\
Flood1 & 1 & 266 & -- & 0.0181 & \textbf{0.0164} & 0.0164 & 0.0181 & 0.0165 & 0.0165 \\
Flood2 & 1 & 266 & -- & 0.0185 & 0.0176 & 0.0176 & 0.0185 & 0.0176 & \textbf{0.0176} \\
Flood3 & 1 & 266 & -- & 0.0214 & \textbf{0.00946} & 0.00947 & 0.0214 & 0.00963 & 0.00963 \\
GasSensorAc & 1 & 7500 & -- & 0.268 & 0.25 & 0.249 & 0.271 & \textbf{0.246} & 0.246 \\
GasSensorEth & 1 & 7500 & -- & 0.299 & 0.29 & 0.292 & 0.301 & \textbf{0.284} & 0.292 \\
MetroTraffic & 4 & 24 & \textbf{535} & 535 & 535 & 535 & 535 & 535 & 535 \\
SierraSnow & 3 & 30 & 35.4 & 33.5 & 24.8 & 25 & 33.7 & \textbf{24.3} & 24.5 \\
\bottomrule
\end{tabular}
\end{adjustbox}
\caption{Optimal RMSE across $\lambda$ of SKRR applied to the TSER dataset for different signature kernel methods. Displayed dimension excludes added time channel.}
\label{table:RMSE-kernel}
\end{table}
\paragraph{AI usage statement.} Generative AI tools, specifically Claude (Opus 5.5 and earlier versions), were used during the preparation of this manuscript to assist with language
editing, literature searches, numerical experiments, and feedback on mathematical arguments. The authors take full responsibility for the correctness and content of the manuscript.
\paragraph{Acknowledgments.}
Both authors thank Davit Gogolashvili for helpful discussions.
CB gratefully acknowledges support by the Deutsche Forschungsgemeinschaft
(DFG) CRC/TRR 388 “Rough Analysis, Stochastic Dynamics and Related Fields”, Projects B01 and B03.
JK gratefully acknowledges financial support by the DFG through the IRTG 2544 “Stochastic Analysis in Interaction”.
\newpage
\appendix
\section{Appendix} \label{Appendix}
The appendix is organized as follows. 
In Section \ref{sec:proofs}, we discuss the proofs of the results in the main body.
Then, in Section \ref{section_sub_mercer}, we discuss Mercer representations 
for the signature kernel, which are closely related to the theory of the integral operator $\mc{T}$, but not directly relevant to the finite-sample analysis.
\subsection{Proofs}\label{sec:proofs}
\subsubsection{Computing the approximation error}\label{sec:approx-error:proofs}
In this section, we discuss the proof of Theorem \ref{theorem_approximation_error} and Corollary \ref{corollary:approx-error:Brownian}, dealing with the approximation error.  The proof of the former is based on a spectral decomposition of the operator $\mc{T}$, whereas the latter relies on the observation that iterated Itô integrals of Brownian motion are orthogonal in $L^2$.

Recall that we introduced $\mc{T}$ as an operator from $L^2 (\nu_{\mc{X}})\to \mc{H}$. We define $\mc{T}_\mch:\mch\to \mch$ analogously for $f\in\mc{H}$ by
\begin{align}\label{eq_def_T_H}
    \mc{T}_{\mc{H}} f \coloneqq 
    \int_\mcx k(\cdot,\bx') f(\bx') \,\nu_\mcx(\mathrm{d}\bx'),
\end{align}
which we observe is positive, self-adjoint, and compact, provided the integrability condition in Assumption \ref{ass_kernel_diagonal} holds.
Further, if we want to distinguish between elements of $\mc{H}$ and their $L^2 (\nu_{\mc{X}})$-equivalence classes, given $f\in \mc{H}$, we write $[f]_{\sim}$ for its $L^2 (\nu_{\mc{X}})$-equivalence class.

The next result is a direct consequence of the spectral Theorem for compact, self-adjoint, positive operators and is central to proving Theorem \ref{theorem_approximation_error}. It can be derived as a consequence of \cite[Lemma 2.12]{SteinwartScovel2012}.
\begin{lemma}\label{lemma:steinwart-scovel}
Consider the setup of Lemma \ref{lemma_RKHS} and $\nu_\mcx\in \mc{P}(\mcx)$ fulfilling the integrability condition in Assumption \ref{ass_kernel_diagonal}, so that $\mc{H}$ is compactly embedded into $L^2 (\nu_{\mc{X}})$. 

Then, there exist an at most 
countable family $(\mu_i)_{i\in I}\subset (0,\infty)$ converging to $0$ with $\mu_1\geq \mu_2 \geq \dots>0$ and a family of functions $(e_i)_{i\in I}\subset \mc{H}$ such that
\begin{itemize}
    \item $(\sqrt{\mu_i}e_i)_{i\in I}$ is an ONS in $\mch$ of eigenvectors diagonalizing $\mc{T}_\mch$ with eigenvalues $(\mu_i)_{i\in I}$ and
    \item $([e_i]_{\sim})_{i\in I}$ is an ONS in $L^2 (\nu_{\mc{X}})$ of eigenvectors of $\mc{T}$ with eigenvalues $(\mu_i)_{i\in I}$, in the sense that for $f\in L^2 (\nu_{\mc{X}})$ we have
    \begin{align}\label{eq_T_k_L2_classical_2}
        \mc{T} f = \sum_{i\in I}\mu_i \langle f,[e_i]_{\sim}\rangle_{L^2 (\nu_{\mc{X}})} e_i.
    \end{align}
\end{itemize}
\end{lemma}
\begin{remark}\label{remark_embedding_kernel}
    Observe that in the case of the signature kernel, for $\mch$ dense in $L^2(\nu_\mcx)$ with $\textnormal{dim}(L^2(\nu_\mcx))=\infty$, we have $I=\N$ and the set $([e_i]_\sim)_{i\in \N}$ is not just an ONS but an ONB of $L^2(\nu_\mcx)$, since $\textnormal{ker}(\mc{T})=\{0\}$. Indeed, for $f\in L^2(\nu_\mcx) \setminus \{0\}$ we have
    \begin{align}\label{eq_kernel_T}
        \langle f,[\mc{T} f]_{\sim}\rangle_{L^2 (\nu_{\mc{X}})} 
        =\int_{\mc{X}}\int_{\mc{X}}f(\bx)f(\bx')\,k(\bx,\bx')\,\nu_{\mc{X}}(\mathrm{d}\bx')\,\nu_{\mc{X}}(\mathrm{d}\bx).
    \end{align}
    Since the above equals zero only if
    we have $\int_{\mc{X}} f(\bx) \sig(\bx)^w\,\nu_\mcx(\mathrm{d}\bx)=0$ for all words $w$ (namely multi-indices), we find that $f\notin \ker (\mc{T})$ if $\mch$ is dense in $L^2(\nu_\mcx)$.
\end{remark}

We now prove Theorem \ref{theorem_approximation_error}.
\begin{proof}[Proof of Theorem \ref{theorem_approximation_error}]
    Denote by $\{e_i\}_{i\in \N} \subset \mch$ the family of eigenvectors associated to $\mc{T}$, respectively $\mc{T}_\mch$, provided by Lemma \ref{lemma:steinwart-scovel}.
    Taking into account the discussion in Remark \ref{remark_embedding_kernel}, $\{[e_i]_\sim\}_{i\in \N}$ is an ONB of $L^2(\nu_\mcx)$.
    As a consequence of the representation $f^{\lambda} = (\mc{T}+\lambda)^{-1}\mc{T} f^{\ast}$, see also \cite[Proposition 3.3]{BayerBelomestnyButkovskyEtAl2024}, one thus has
    \begin{align}\notag
        f^{\lambda} -f^{\ast}= ((\mc{T}+\lambda)^{-1}\mc{T} -\textnormal{Id})f^{\ast}
        = -\sum_{i\in \mb{N}} \frac{\lambda}{\mu_i+\lambda} \langle [e_i]_\sim,f^{\ast}\rangle_{L^2(\nu_{\mc{X}})}[e_i]_\sim,
    \end{align}
    which implies \eqref{eq_L2_error_estimate}.
    One obtains \eqref{eq_H_error_estimate} in a similar fashion.
\end{proof}

We now derive Corollary \ref{corollary:approx-error:Brownian} as a consequence of Theorem \ref{theorem_approximation_error}.
\begin{proof}[Proof of Corollary \ref{corollary:approx-error:Brownian}]
    Note first that any element $f\in \mch$ can be written as $f(\bx) = \langle \ell_f,\sig(\bx)\rangle$ for some $\ell_f\in T((\mb{R}^d))$. 
    Observe that different iterated integrals are orthogonal in $L^2(\nu_\mcx)$.
    Since the $n$-th level iterated Itô integrals have variance $T^n/n!$, we have
    \begin{align}
        \mc{T}_\mch f (\bx) =& \sum_{n\in \mb{N}} \sum_{w\in W_n} \int_\mcx (\sig(\bx')^w)^2 \,\nu_{\mc{X}}(\mathrm{d}\bx')\,\ell_f^w\sig(\bx)^w
        \\ \notag = & \sum_{n\in \mb{N}} \sum_{w\in W_n} \frac{T^n}{n!}\,\ell_f^w\, \sig(\bx)^w.
    \end{align}
    Hence, the eigenfunctions of $\mc{T}_\mch$ with positive eigenvalues are given by the iterated integrals, i.e., for $w\in W_n$, the function $\bx\mapsto \sig(\bx)^w\in \mch$ is an eigenvector of $\mc{T}_\mch$ with eigenvalue $\frac{T^n}{n!}$ and the collection $\{\sig^w\}_{w \in W}$ is an eigenbasis of $\textnormal{ker}(\mc{T})^\perp=\bar{\mch}$. Here, we treat the map $\bx \mapsto \sig (\bx)^{\emptyset}=1$ as the zero-order iterated integral.
    Using that 
    \begin{equation*}
        \big\|f^{\ast}-f^{\lambda}\big\|_{L^2(\nu_{\mc{X}})}^2
        =
        \big\|f^{\ast}-P_{\bar{\mch}} f^{\ast}\big\|_{L^2(\nu_{\mc{X}})}^2
        +
        \big\|P_{\bar{\mch}} f^{\ast}-f^{\lambda}\big\|_{L^2(\nu_{\mc{X}})}^2,
    \end{equation*}
    the Corollary follows as a direct consequence of Lemma \ref{lemma:steinwart-scovel} and the proof of Theorem \ref{theorem_approximation_error}, switching $f^{\ast}$ for $P_{\bar{\mch}} f^{\ast}$.
\end{proof}

\subsubsection{Computing the finite-sample error under a Bernstein condition}\label{section_sub_finite_sample_error}
The objective of this section is to prove Theorems \ref{theorem:excess-risk:Bernstein} and \ref{theorem:rates:Bernstein}, dealing with the finite-sample error under Assumption \ref{assumption:finitesample:bernstein}. As the latter relies on Theorem \ref{theorem:excess-risk:Bernstein}, we start establishing this one. 

We further rely on the following classical representation of the error criterion $\mc{E}_{\textnormal{rel}}$ used in Theorem \ref{theorem:excess-risk:Bernstein} 
in terms of an $\mc{H}$-norm, see for example \cite[Proposition 1]{Caponnetto2007OptimalRates}.
\begin{lemma}\label{lemma_E_error_rep}
    Provided $\nu_{\mb{R}}\in \mc{P}_2(\R)$, for any $f\in L^2(\nu_{\mc{X}})$ we have
    \begin{align}\label{eq_E_error_rep}
        \mc{E}[f]-\mc{E}[f^{\ast}] = 
        \left\|
        \sqrt{\mc{T}}(f-f^{\ast})
        \right\|_{\mch}^2.
    \end{align}
\end{lemma}
We now prove Theorem \ref{theorem:excess-risk:Bernstein}.
\begin{proof}[Proof of Theorem \ref{theorem:excess-risk:Bernstein}]
    We follow the strategy of the proof of Theorem 4 from \cite{Caponnetto2007OptimalRates}, making suitable adjustments to deal with the unboundedness of the signature kernel. 
    We split the proof into four steps, as is done in \cite{Caponnetto2007OptimalRates}.
    
    In the following, 
    we denote a collection of $N$ i.i.d. random variables with law $\nu$ by $(Z_i)_{1\leq i \leq N}=(\bX_i,Y_i)_{1\leq i \leq N}$.
    We will denote the collections of samples by the index $1:N$, e.g. we will write $Z_{1:N} = (Z_i)_{1\leq i \leq N}$, $\bX_{1:N} = (\bX_i)_{1\leq i \leq N}$, and $Y_{1:N} = (Y_i)_{1\leq i \leq N}$.
    
    For $\bX\sim \nu_\mcx$, we define $\mc{T}_{\bX}$ to be the random operator that maps $f\in \mch$ to
    \begin{align}
        \mc{T}_{\bX}f \coloneqq f(\bX) \, k_{\bX}= f(\bX) \, k(\bX,\cdot)\in \mch.
    \end{align}
    Given $f\in L^2(\nu_{\mc{X}})$, we use the same notation for the $\mch$-valued random variable $\mc{T}_{\bX}f = f(\bX) \, k(\bX,\cdot)$, which does not depend on the choice of representative of $f$ up to almost sure equality.
    
    If we consider sample averages of operators depending on the single samples, we will also indicate this by the index $1:N$, e.g. we write
    \begin{align}
        \mc{T}_{\bX_{1:N}} \coloneqq \frac{1}{N} \sum_{1\leq i \leq N} \mc{T}_{\bX_i}.
    \end{align}
    We will also rely in the following on the fact that
    \begin{align}
        f^{\lambda}_N = \left\langle (k(\bX_i,\cdot))_{1\leq i \leq N},(K+\lambda N)^{-1} Y_{1:N} \right\rangle_{\mb{R}^N},
    \end{align}
    where here $K= (k(\bX_i,\bX_j))_{1\leq i,j\leq N}$,
    can be expressed equivalently in terms of the operator $\mc{T}_{\bX_{1:N}}$ by
    \begin{align}
        f^{\lambda}_N= \left(\mc{T}_{\bX_{1:N}} + \lambda\right)^{-1}\frac{1}{N}\sum_{i=1}^N k_{\bX_i}Y_i.
    \end{align}
    We now begin proving the Theorem.
    
    \textbf{Step 1} As shown in Lemma \ref{lemma_E_error_rep}, it holds that 
    \begin{align}\label{eq_error_rep}
        \mc{E}_{\textnormal{rel}}[f_{N}^\lambda]
        =\mc{E}[f_{N}^\lambda]-\mc{E}[f^{\ast}]
        = \| \sqrt{\mc{T}} (f_{N}^\lambda-f^{\ast})\|_{\mch}^2.
    \end{align}
    Writing $g_{Z_{1:N}} = \frac{1}{N}\sum_{i=1}^N k_{\bX_i}Y_i$ and $\mc{T}_{\bX_{1:N}} f = \frac{1}{N} \sum_{1\leq i \leq N} \mc{T}_{\bX_i}f$, we define the quantities
    \begin{align}\label{eq:def:S}
        &S_1(\lambda, Z_{1:N})
        =
        \left\|
        \sqrt{\mc{T}}\,(\mc{T}_{\bX_{1:N}} + \lambda)^{-1}\big(g_{Z_{1:N}} - \mc{T}_{\bX_{1:N}} f^{\ast}\big)
        \right\|_{\mch}^{2},
        \\ \notag &
        S_2(\lambda, Z_{1:N})
        =
        \left\|
        \sqrt{\mc{T}}\,(\mc{T}_{\bX_{1:N}} + \lambda)^{-1}(\mc{T} - \mc{T}_{\bX_{1:N}})\big(f^\lambda - f^{\ast}\big)
        \right\|_{\mch}^{2}.
    \end{align}
    Further recall that the residual was defined as
    \[
        \mc{A} (\lambda) =\mc{E}[f^\lambda]-\mc{E}[f^{\ast}]=\|\sqrt{\mc{T}}(f^\lambda-f^{\ast})\|_{\mch}^2.
    \]
    Then, using $(a+b+c)^2\leq 3(a^2+b^2+c^2)$ and the identity \eqref{eq_error_rep},
    we obtain
    \begin{align}\label{eq:Etrianglebound}
        \mc{E}[f_{N}^\lambda]-\mc{E}[f^{\ast}]
        \le 3\big( \mc{A}(\lambda) + S_1(\lambda, Z_{1:N}) + S_2(\lambda, Z_{1:N}) \big).
    \end{align}
    In the next two steps, we establish bounds on the quantities $S_1(\lambda, Z_{1:N}),S_2(\lambda, Z_{1:N})$.
    
    \textbf{Step 2:} To bound $S_2$, note that by definition of the linear operator norm, one has
    \begin{align}\label{eq_S2_first_bound}
        S_2(\lambda, Z_{1:N})
        \leq 
        \big\|
        \sqrt{\mc{T}}\,(\mc{T}_{\bX_{1:N}} + \lambda)^{-1}\big\|_{\mc{L}(\mch)}^{2}
        \cdot \big\|(\mc{T} - \mc{T}_{\bX_{1:N}})\big(f^\lambda - f^{\ast}\big)
        \big\|_{\mch}^{2}.
    \end{align}
    We now bound the two factors on the right hand side of \eqref{eq_S2_first_bound}.
    
    \textbf{Step 2.1:} We start bounding the first term on the right hand side of \eqref{eq_S2_first_bound}. By a Neumann series expansion argument, it follows that the inequality
    \begin{align}\label{eq_bound_theta}
        \Theta(\lambda, Z_{1:N}) &
        \coloneqq
        \left\|(\mc{T} + \lambda)^{-1}(\mc{T} - \mc{T}_{\bX_{1:N}})\right\|_{\mathcal{L}(\mch)}
        \le 1/2,
    \end{align}
    implies
    \begin{align}
        \big\|
        \sqrt{\mc{T}}\,(\mc{T}_{\bX_{1:N}} + \lambda)^{-1}\big\|_{\mc{L}(\mch)}
        \leq \frac{1}{\sqrt{\lambda}}.
    \end{align}
    Indeed, since $(\mc{T} + \lambda)^{-1}$ and $\mc{T} - \mc{T}_{\bX_{1:N}}$ are self-adjoint, 
    we have
    \begin{equation*}
        \Theta(\lambda, Z_{1:N}) 
        =
        \left\|(\mc{T} - \mc{T}_{\bX_{1:N}}) (\mc{T} + \lambda)^{-1}\right\|_{\mathcal{L}(\mch)},
    \end{equation*}
    so 
    \begin{align*}
        \big\|
        \sqrt{\mc{T}}\,(\mc{T}_{\bX_{1:N}} + \lambda)^{-1}\big\|_{\mc{L}(\mch)}
        = &
        \big\|
        \sqrt{\mc{T}}\,(\mc{T} + \lambda)^{-1} (I-(\mc{T}-\mc{T}_{\bX_{1:N}})(\mc{T}+\lambda)^{-1})^{-1}\big\|_{\mc{L}(\mch)}
        \\ \notag = & 
        \big\|
        \sqrt{\mc{T}}\,(\mc{T} + \lambda)^{-1} \sum_{n\in \N}\left((\mc{T}-\mc{T}_{\bX_{1:N}})(\mc{T}+\lambda)^{-1}\right)^n\big\|_{\mc{L}(\mch)}
        \\ \notag \leq & 
        \big\|
        \sqrt{\mc{T}}\,(\mc{T} + \lambda)^{-1} 
        \big\|_{\mc{L}(\mch)}
        \sum_{n\in \N}\big\|
        (\mc{T}-\mc{T}_{\bX_{1:N}})(\mc{T}+\lambda)^{-1}\big\|_{\mc{L}(\mch)}^n
        \\ \notag \leq &  2 \big\|
        \sqrt{\mc{T}}\,(\mc{T} + \lambda)^{-1} 
        \big\|_{\mc{L}(\mch)}.
    \end{align*}
    Since the spectral Theorem ensures $\big\|
    \sqrt{\mc{T}}\,(\mc{T} + \lambda)^{-1} 
    \big\|_{\mc{L}(\mch)} \leq \frac{1}{2\sqrt{\lambda}}$, the implication follows.
    
    Relying on Theorem \ref{theorem:Bernstein}, we thus establish next that for $N$ and $\lambda$ fulfilling \eqref{eq_condition_params}, the inequality \eqref{eq_bound_theta} holds. 
    We do so by applying Theorem \ref{theorem:Bernstein} to the operator valued random variable
    \begin{align}
        \zeta_1 = (\mc{T}+\lambda)^{-1} \mc{T}_{\bX},
    \end{align}
    where $\bX\sim \nu_{\mc{X}}$, so
    \begin{align}
        \mb{E}[\zeta_1 ]= (\mc{T} + \lambda)^{-1}\mc{T}
        \,\,\,\,\,\,\textnormal{and}\,\,\,\,\,\,\frac{1}{N}\sum_{i=1}^N \zeta_1^i = 
        (\mc{T}+\lambda)^{-1} \mc{T}_{\bX_{1:N}},
    \end{align}
    where $\zeta^i_1$ are $N$ i.i.d. copies of $\zeta_1$ obtained by replacing $\bX$ by $\bX_i\sim \nu_{\mcx}$.
    To verify the moment condition of Theorem \ref{theorem:Bernstein}, we compute the moments of the Hilbert-Schmidt norm, which upper bounds the linear operator norm. To this end, note that
    \begin{align}\notag
        & \mb{E}\big[\|\zeta_1-\mb{E}[\zeta_1]\|^m_{\HS}\big]
        =
        \mb{E}\big[\|(\mc{T}+\lambda)^{-1}(\mc{T}_{\bX}-\mc{T})\|^m_{\HS}\big]
        \\ \notag \leq &
        \mb{E}\big[\|(\mc{T}+\lambda)^{-1}\|^m_{\mc{L}(\mch)}
        \cdot \|\mc{T}_{\bX}-\mc{T}\|^m_{\HS}\big]
        \\ \notag \leq &
        \frac{1}{\lambda^m}\mb{E}\big[
        \|\mc{T}_{\bX}-\mc{T}\|^m_{\HS}\big]
        \\ \notag \leq &
        \frac{1}{\lambda^m}
        \big( \mb{E}\big[
        \|\mc{T}_{\bX}\|^m_{\HS}\big]^{1/m}+
        \|\mc{T}\|_{\HS}
        \big)^m
        \\ \label{eq_step_21} \leq &
        \frac{2^m}{\lambda^m}
        \mb{E}\big[
        \|\mc{T}_{\bX}\|^m_{\HS}\big],
    \end{align}
    where in the last inequality we used that by Jensen $\|\mc{T}\|_{\HS} \leq \mb{E}\big[
        \|\mc{T}_{\bX}\|^m_{\HS}\big]^{1/m}$. 
    We now compute $\mb{E}\big[
        \|\mc{T}_{\bX}\|^m_{\HS}\big]$. 
        By Riesz and the reproducing property, we have
        $\mc{T}_\bX(f) = \langle k_\bX,f \rangle_{\mch} k_\bX$,
        so we find
        \begin{align}\label{eq_T_x_computation}
            \|\mc{T}_\bX\|_{\HS}= \|k_\bX\|^2_\mch = k(\bX,\bX).
        \end{align}
        In particular, we thus have by \eqref{eq_step_21} that
        \begin{align}\label{eq:zeta1moments}
            \mb{E}\big[\big\|\zeta_1-\mb{E}[\zeta_1]\big\|^m_{\HS}\big]
            \leq
            \frac{2^m}{\lambda^m}
            \mb{E}\big[
            k(\bX,\bX)^m\big].
        \end{align}
        By Assumption \ref{assumption:finitesample:bernstein}, the above implies that for all $m\geq 2$ we have
    \begin{align}\notag
        \mb{E}\big[\big\|\zeta_1-\mb{E}[\zeta_1]\big\|^m_{\HS}\big]
        \leq 
        m! \big(\frac{2C_{\nu}}{\lambda}\big)^m
        \leq \frac{1}{2}
        m! 
        \underbrace{\big(\frac{3C_{\nu}}{\lambda}\big)^2}_{\sigma^2}
        \underbrace{\big(\frac{2C_{\nu}}{\lambda}\big)^{m-2}}_{L^{m-2}}.
    \end{align}
    Hence, Theorem \ref{theorem:Bernstein} implies that with probability at least $1-\eta/3$ we have
    \begin{align}\notag
        \|(\mc{T}+\lambda)^{-1}(\mc{T}_{\bX_{1:N}}-\mc{T})\big\|_{\mc{L}(\mch)}
        \leq \|(\mc{T}+\lambda)^{-1}(\mc{T}_{\bX_{1:N}}-\mc{T})\big\|_{\HS}
        \leq \frac{4C_\nu}{\lambda N} \log \frac{6}{\eta}
        +
        \frac{3C_{\nu}}{\lambda \sqrt{N}} 
        \sqrt{2 \log \frac{6}{\eta}}
        .
    \end{align}
    For $N$ and $\lambda$ fulfilling \eqref{eq_condition_params}, the above right hand side is bounded by $1/2$, so \eqref{eq_bound_theta} holds. 
    
    \textbf{Step 2.2.} We bound $\|(\mc{T}-\mc{T}_{\bX_{1:N}})(f^\lambda-f^{\ast})\|_{\mch}$ by applying Theorem \ref{theorem:Bernstein} to the $\mch$ valued r.v. $\zeta_2$ defined by
    \begin{align}
        \zeta_2 =\mc{T}_{\bX} (f^\lambda-f^{\ast}),
    \end{align}
    where again $\bX\sim \nu_{\mcx}$.
    Then, we have
    \begin{align}
        \mb{E} \left[\zeta_2\right]=\mc{T}(f^\lambda-f^{\ast})\,\,\,\,\,\,\textnormal{and}\,\,\,\,\,\,\frac{1}{N}\sum_{i=1}^N \zeta_2^i = \mc{T}_{\bX_{1:N}}(f^\lambda-f^{\ast}),
    \end{align}
    where $\zeta^i_2$ are $N$ i.i.d. copies of $\zeta_2$ obtained by replacing $\bX$ by $\bX_i\sim \nu_{\mcx}$.
    To bound the associated moments, note that
    \begin{align}\notag
        \mb{E}\big[\|\zeta_2-\mb{E}[\zeta_2]\|^m_{\mch}\big]
        = &
        \mb{E}\big[\|(\mc{T}_{\bX}-\mc{T}) (f^\lambda-f^{\ast})\|^m_{\mch}\big]
        \\ \notag \leq &
        \mb{E}\big[
        \|\mc{T}_{\bX}-\mc{T}\|_{\mc{L}(\mch)}^m \|f^\lambda-f^{\ast}\|_{\mch}^m
        \big]
        \\ \notag
        =
        &\mb{E}\big[
        \|\mc{T}_{\bX}-\mc{T}\|_{\mc{L}(\mch)}^m 
        \big]
        \|f^\lambda-f^{\ast}\|_{\mch}^m
        \\ \notag
        = 
        &\mb{E}\big[
        \|\mc{T}_{\bX}-\mc{T}\|_{\mc{L}(\mch)}^m 
        \big]
        {\mc{B}(\lambda)}^{m/2}.
    \end{align}
    Analogously to the moment bound for $\zeta_1$, we thus find by Assumption \ref{assumption:finitesample:bernstein} that for all $m\geq 1$ we have
    \begin{align}
        \mb{E}\big[\|\zeta_2-\mb{E}[\zeta_2]\big\|^m_{\mch}\big]
        \leq
        \frac{1}{2}
        m! 
        \big(3C_\nu
        {\mc{B}(\lambda)}^{1/2}
        \big)^2
        \big(2C_\nu {\mc{B}(\lambda)}^{1/2}\big)^{m-2}.
    \end{align}
    Hence, due to Theorem \ref{theorem:Bernstein}, with probability at least $1-\eta/3$, we have
    \begin{align}\label{eq:S2:boundpart2}
        \| (\mc{T}-\mc{T}_{\bX_{1:N}})(f^\lambda-f^{\ast}) \|_{\mch}
        \leq  
        \frac{4C_\nu {\mc{B}(\lambda)}^{1/2}}{N}
        \log \frac{6}{\eta}
        +\frac{3C_\nu
        {\mc{B}(\lambda)}^{1/2}}{\sqrt{N}}
        \sqrt{2\log \frac{6}{\eta}}.
    \end{align}
    \textbf{Step 2.3:} Combining the bounds on the quantities in \eqref{eq_S2_first_bound} derived in steps 2.1 and 2.2, we find that in total, provided $\lambda,N$ fulfill \eqref{eq_condition_params}, with probability at least $1-2\eta/3$ we have
    \begin{align}\notag
        S_2(\lambda, Z_{1:N})
        \leq &
        \big\|
        \sqrt{\mc{T}}\,(\mc{T}_{\bX_{1:N}} + \lambda)^{-1}\big\|_{\mc{L}(\mch)}^{2}
        \cdot \big\|(\mc{T} - \mc{T}_{\bX_{1:N}})\big(f^\lambda - f^{\ast}\big)
        \big\|_{\mch}^{2}
        \\ \notag \leq &
        \frac{1}{\lambda}
        \bigg(\frac{4C_\nu {\mc{B}(\lambda)}^{1/2}}{N}
        \log \frac{6}{\eta}
        +\frac{3C_\nu
        {\mc{B}(\lambda)}^{1/2}}{\sqrt{N}}
        \sqrt{2\log \frac{6}{\eta}}
        \bigg)^2
        \\ \notag \leq &
        \frac{2}{\lambda}
        \log \frac{6}{\eta}
        \Big(\frac{4^2 C_{\nu}^2 {\mc{B}(\lambda)}}{N^2}
        \log \frac{6}{\eta}
        +\frac{2 \cdot 3^2C_{\nu}^2
        {\mc{B}(\lambda)}}{N}\Big)
        \\ \notag \leq &
        32
        \frac{C_{\nu}^2 {\mc{B}(\lambda)}}{\lambda N^2}
        \log^2 \frac{6}{\eta}
        +36 \frac{C_{\nu}^2
        {\mc{B}(\lambda)}}{\lambda N}
        \log \frac{6}{\eta}
        .
    \end{align}
    \textbf{Step 3:} We now bound $S_1(\lambda,Z_{1:N})$. Note first that
    \begin{align}\label{eq_bound_S1}
        & S_1(\lambda, Z_{1:N})
        =
        \left\|
        \sqrt{\mc{T}}\,(\mc{T}_{\bX_{1:N}} + \lambda)^{-1}\big(g_{Z_{1:N}} - \mc{T}_{\bX_{1:N}} f^{\ast}\big)
        \right\|_{\mch}^{2}
        \\ \notag \leq &
        \left\|
        \sqrt{\mc{T}}\,(\mc{T}_{\bX_{1:N}} + \lambda)^{-1}(\mc{T}+\lambda)^{1/2}\right\|_{\mc{L}(\mch)}^2\left\|
        (\mc{T}+\lambda)^{-1/2}
        \big(g_{Z_{1:N}} - \mc{T}_{\bX_{1:N}} f^{\ast}\big)
        \right\|_{\mch}^{2}.
    \end{align}
    We now bound the two terms on the right hand side of the above inequality using Theorem \ref{theorem:Bernstein}.
    
    \textbf{Step 3.1} For the first term in \eqref{eq_bound_S1}, by a Neumann series expansion argument, similarly as done in step 2.1, it is shown in \cite{Caponnetto2007OptimalRates} that under the conditions imposed in step 2, which hold with probability at least $1-2\eta /3$ provided $\lambda,N$ fulfill \eqref{eq_condition_params}, we have that 
    \begin{align}
        \left\|
        \sqrt{\mc{T}}\,(\mc{T}_{\bX_{1:N}} + \lambda)^{-1}(\mc{T}+\lambda)^{1/2}\right\|_{\mc{L}(\mch)}^2 \leq 4.
    \end{align}
    In particular, the above thus holds with probability at least $1-2\eta /3$.
    
    \textbf{Step 3.2} We now bound the second contribution in \eqref{eq_bound_S1} using concentration arguments via Theorem \ref{theorem:Bernstein}.
    Introduce the $\mch$-valued r.v. $\zeta_3$ by
    \begin{align}
        \zeta_3
        = (\mc{T}+\lambda)^{-1/2}k_{\bX} (Y-f^{\ast}(\bX)),
    \end{align}
    where $(\bX,Y)\sim \nu$. Note that
    \begin{align}
        \mb{E}[\zeta_3]=0
        \,\,\,\,\,\,\textnormal{and}\,\,\,\,\,\,\frac{1}{N}\sum_{i=1}^N \zeta_3^i = 
        (\mc{T}+\lambda)^{-1/2}
        \big(g_{Z_{1:N}} - \mc{T}_{\bX_{1:N}} f^{\ast}\big),
    \end{align}
    where $\zeta_3^i$ are $N$ i.i.d. copies of $\zeta_3$ obtained by replacing $(\bX,Y)$ by $Z_i=(\bX_i,Y_i)\sim \nu$.
    Recalling that Assumption \ref{assumption:finitesample:bernstein} states that for all $m\geq 2$ we have
    \begin{align}\label{eq_ass_nu_recall}
       \int_{\mb{R}}|y-f^{\ast} (\bx)|^m\,\nu(\mathrm{d}y|\bx) \leq \frac{1}{2} (m/2)!\,\Sigma^2M^{m-2},
    \end{align}
    we can control the moments of $\zeta_3$ by
    \begin{align}\notag
        \mathbb{E}\big[\|\zeta_3\|_{\mch}^m\big]
        &=
        \int_{\mcx\times \R}
        \left(
        \left\langle
        (\mc{T}+\lambda)^{-1}k_\bx \big(y - f^{\ast}(\bx)\big),
        k_\bx\, (y - f^{\ast}(\bx))
        \right\rangle_\mch
        \right)^{m/2}
        \, \nu(\mathrm{d}\bx,\mathrm{d}y)
        \\ \notag
        &\le
        \int_{\mcx}
        \|(\mc{T}+\lambda)^{-1}\|_{\mathcal{L}(\mch)}^{m/2}
        \|k_\bx\|_\mch^{m}
        \left(
        \int_{\R}
        |y - f^{\ast}(\bx)|^{m}
        \, \nu(\mathrm{d}y|\bx)
        \right)
        \, \nu_{\mc{X}}(\mathrm{d}\bx)
        \\ \label{eq:zeta3bound}
        &=\frac{1}{\lambda^{m/2}}
        \int_{\mcx}
        k(\bx,\bx)^{m/2}
        \left(
        \int_{\R}
        |y - f^{\ast}(\bx)|^{m}
        \, \nu(\mathrm{d}y|\bx)
        \right)
        \, \nu_{\mc{X}}(\mathrm{d}\bx).
    \end{align}
    Plugging \eqref{eq_ass_nu_recall} in the above, we thus obtain for even $m$ that
    \begin{align}\notag
        \mathbb{E}\big[\|\zeta_3\|_{\mch}^m\big]
        \leq &
        ( m/2)! (\frac{C_{\nu}}{\lambda})^{m/2}
        \frac{1}{2} (m/2)!\,\Sigma^2M^{m-2}
        \\ \notag \leq &
        \frac{1}{2} m! (\frac{C_{\nu}}{\lambda})^{m/2}
        \,\Sigma^2M^{m-2}
        \\ \notag = &
        \frac{1}{2} m! 
        \,\left(\Sigma \sqrt{\frac{C_{\nu}}{\lambda}}\right)^2
        \left(M\sqrt{\frac{C_{\nu}}{\lambda}}\right)^{m-2}
        .
    \end{align}
    For odd $m$, the above implies by Cauchy--Schwarz that
    \begin{align} \notag
        \mathbb{E}\big[\|\zeta_3\|_{\mch}^m\big]
        \leq    &
        \sqrt{\mathbb{E}\big[\|\zeta_3\|_{\mch}^{m+1}\big]\mathbb{E}\big[\|\zeta_3\|_{\mch}^{m-1}\big]}
        \\ \notag \leq &
        \frac{1}{2}m!
        \,\left(\Sigma \sqrt{\frac{2C_{\nu}}{\lambda}}\right)^2
        \left(M\sqrt{\frac{C_{\nu}}{\lambda}}\right)^{m-2}
        .
    \end{align}
    Therefore, by Theorem \ref{theorem:Bernstein}, with probability at least $1-\eta/3$ we obtain 
    \begin{align}
        &
        \left\|
        (\mc{T}+\lambda)^{-1/2}
        \big(g_{Z_{1:N}} - \mc{T}_{\bX_{1:N}} f^{\ast}\big)
        \right\|_{\mch}
        \\ \notag \leq &
        2M\frac{1}{N}\sqrt{\frac{C_\nu}{\lambda}}\cdot\log\frac{6}{\eta} + \Sigma\sqrt{\frac{2C_\nu}{\lambda N}}\cdot\sqrt{2\log\frac{6}{\eta}}
    \end{align}
    Hence, by \eqref{eq_bound_S1} we find
    that, provided $\lambda,N$ fulfill \eqref{eq_condition_params}, with probability at least $1-\eta$ we have 
    \begin{align}
        S_1(\lambda, Z_{1:N}) 
        \leq 32 \Sigma^2 C_\nu\cdot\log \frac{6}{\eta} \frac{1}{\lambda N} + 32 M^2 C_\nu \cdot\log^2 \frac{6}{\eta} \frac{1}{\lambda N^2}.
    \end{align}
    \textbf{Step 4.} Combining the results of the previous steps, for $\lambda,N$ as proposed in \eqref{eq_condition_params}, we obtain that with probability at least $1-\eta$ we have
    \begin{align}
        \mc{E}[f_{N}^\lambda]-\mc{E}[f^{\ast}]
        \le & 
        3\mc{A}(\lambda) + \left(C_\nu \log(6/\eta)\frac{1}{\lambda N}\right)\cdot(96 \Sigma^2 + 108 C_\nu \mc{B}(\lambda)) 
        \\ & + C_\nu \log^2(6/\eta) \frac{1}{\lambda N^2}\cdot(96M^2 + 96 C_\nu \mc{B}(\lambda))
    \end{align}
    which is the claim of the Theorem.
    \end{proof}
    We now use Theorem \ref{theorem:excess-risk:Bernstein} to establish Theorem \ref{theorem:rates:Bernstein}.
    The next classical result contains a control on the terms $\mc{A}(\lambda)$ and $\mc{B}(\lambda)$, appearing in the bound of Theorem \ref{theorem:excess-risk:Bernstein}.
    \begin{lemma}[{\cite[Proposition 3]{Caponnetto2007OptimalRates}}]\label{lemma_A_B_bound}
        In the setup of Theorem \ref{theorem_approximation_error}, we have
        \begin{align}
            \mc{A}(\lambda) = 
            \lambda^2 \sum_{i\in \mb{N}}\langle [e_i]_{\sim},f^{\ast}\rangle_{L^2(\nu_{\mc{X}})}^2 \frac{1}{(\mu_i+\lambda)^2}
        \end{align}
        and, provided $f^{\ast}\in \mch$, we further have
        \begin{align}
            \mc{B}(\lambda)=
            \lambda^2 \sum_{i\in \mb{N}}\langle [e_i]_{\sim},f^{\ast}\rangle_{L^2(\nu_{\mc{X}})}^2 \frac{1}{\mu_i (\mu_i+\lambda)^2}.
        \end{align}
        In particular, we have $\lim_{\lambda\to 0} \mc{A}(\lambda)=0$ and, provided $f^{\ast}\in \mch$, also $\lim_{\lambda\to 0} \mc{B}(\lambda)=0$.
        Further, if $f^{\ast}\in \mch$, the above implies
        \begin{align}\label{eq_bounds_A_B}
            \mc{A}(\lambda) \leq 
            \lambda \|f^{\ast}\|_{\mch}^2
            \,\,\,\,\,\,
            \textnormal{    and   }
            \,\,\,\,\,\,
            \mc{B}(\lambda) \leq \|f^{\ast}\|_{\mch}^2.
        \end{align}
    \end{lemma}
    We now prove Theorem \ref{theorem:rates:Bernstein}.
    \begin{proof}[Proof of Theorem \ref{theorem:rates:Bernstein}]
        Since $f^{\ast}\in \mch$,
        plugging in the bounds \eqref{eq_bounds_A_B} into the bound of Theorem \ref{theorem:excess-risk:Bernstein}, with probability at least $1-\eta$ it holds that
        \begin{align}\notag
            \Erel[f_{N}^{\lambda}]
            \le & 
            3 \lambda\|f^{\ast}\|_{\mch}^2
            + \left(C_\nu \log(6/\eta)\frac{1}{\lambda N}\right)\cdot(96 \Sigma^2 + 108 C_\nu \|f^{\ast}\|_{\mch}^2) 
            \\ & + C_\nu \log^2(6/\eta) \frac{1}{\lambda N^2}\cdot(96M^2 + 96 C_\nu \|f^{\ast}\|_{\mch}^2)
        \end{align}
        provided $\lambda,N$ fulfill the condition \eqref{eq_condition_params}.
        Plugging the proposal for $\lambda_N$ from the Theorem into the above, for which \eqref{eq_condition_params} holds, we obtain the claimed bound on $\Erel[f_{N}^\lambda]$.
    \end{proof}
\subsubsection{Computing the finite-sample error for heavy-tailed signature kernels}\label{sec:heavytailed:proofs}
We now prove the results for the case of heavy-tailed signature kernels, i.e., Theorems \ref{theorem:errorbound:heavytail},
    \ref{theorem:rates:heavytail}, and
    \ref{theorem:consistency:heavytail}.
    Instead of Bernstein, we will use the Fuk-Nagaev inequality \cite{Fuk1973, Nagaev1979}.
    
    We start with the proof of Theorem \ref{theorem:errorbound:heavytail}.
    \begin{proof}[Proof of Theorem \ref{theorem:errorbound:heavytail}]
        The proof follows along the lines of the proof 
        of Theorem \ref{theorem:excess-risk:Bernstein}.
        We start collecting results derived in the proof of Theorem \ref{theorem:excess-risk:Bernstein}. Recall first that in \eqref{eq:Etrianglebound} we established
        \begin{align}\label{eq:Etrianglebound:heavytailed}
            \mc{E}[f_{N}^\lambda]-\mc{E}[f^{\ast}]
            \le 3\big( \mc{A}(\lambda) + S_1(\lambda, Z_{1:N}) + S_2(\lambda, Z_{1:N}) \big),
        \end{align}
        where the residual $\mc{A}(\lambda)$ is defined in \eqref{eq:residual} and $S_1,S_2$ are defined in \eqref{eq:def:S}. 
        
        We derived in Steps 2 and 3 of the proof of Theorem \ref{theorem:excess-risk:Bernstein} that
        \begin{align}\label{eq:S2ineq}
            S_2(\lambda, Z_{1:N})
            \leq &
            \big\|
            \sqrt{\mc{T}}\,(\mc{T}_{\bX_{1:N}} + \lambda)^{-1}\big\|_{\mc{L}(\mch)}^{2}
            \cdot 
            \big\|(\mc{T} - \mc{T}_{\bX_{1:N}})\big(f^\lambda - f^{\ast}\big)
            \big\|_{\mch}^{2},
            \\ \label{eq:S1ineq}
            S_1(\lambda, Z_{1:N})
            \leq &
            \left\|
            \sqrt{\mc{T}}\,(\mc{T}_{\bX_{1:N}} + \lambda)^{-1}(\mc{T}+\lambda)^{1/2}\right\|_{\mc{L}(\mch)}^2
            \left\|
            (\mc{T}+\lambda)^{-1/2}
            \big(g_{Z_{1:N}} - \mc{T}_{\bX_{1:N}} f^{\ast}\big)
            \right\|_{\mch}^{2}.
        \end{align}
        We bound the terms $S_2,S_1$ following the framework of Theorem \ref{theorem:excess-risk:Bernstein}, replacing Bernstein by Fuk-Nagaev.
        
        First, we used that
        the first quantity on the right hand side of \eqref{eq:S2ineq} is bounded by $1/\lambda$
        with probability of at least $1-\eta/3$, provided that
        \begin{align}
            \left\|\frac{1}{N}\sum_{i=1}^N \zeta_1^i - \E[\zeta_1]\right\|_{\mc{L}(\mch)} \leq \frac{1}{2},
        \end{align}
        where $\zeta_1$ is as introduced in the proof of Theorem \ref{theorem:excess-risk:Bernstein}.
        As a consequence of \eqref{eq:zeta1moments} and Theorem \ref{theorem:fuknagaev} applied to the $m$-th moment, we find that for choices of $\lambda,N$ as in condition \eqref{condition:errorbound:heavytail}, the above is indeed bounded by $1/2$ with probability at least $1-\frac{\eta }{3}$. 
        
        The second quantity on the right hand side of \eqref{eq:S2ineq} equals \mbox{$\frac{1}{N}\sum_{i=1}^N \zeta_2^i-\mb{E}[\zeta_2]$,} for $\zeta_2$ introduced as in the proof of Theorem \ref{theorem:excess-risk:Bernstein}. For $\zeta_2$ we showed
        \begin{align}
            \mb{E}\big[\|\zeta_2-\mb{E}[\zeta_2]\|^m_{\mch}\big]
            \leq 
            (2{\mc{B}(\lambda)}^{1/2})^m
            \mb{E}\big[
            k(\bX,\bX)^m\big]
            \leq (2{\mc{B}(\lambda)}^{1/2})^m C_m.
        \end{align}
        Hence, Theorem \ref{theorem:fuknagaev} implies that with probability at least $1-\frac{\eta}{3}$ we have
        \begin{align}\label{eq:S2:boundpart2:heavy} 
            \big\|(\mc{T} - \mc{T}_{\bX_{1:N}})\big(f^\lambda - f^{\ast}\big)
            \big\|_{\mch}
            \leq 
            \max\left\{
            \left(\frac{6c_1}{\eta N^{m-1}}
            (2{\mc{B}(\lambda)}^{1/2})^m C_m
            \right)^{1/m},
            \sqrt{{\mc{B}(\lambda)} C_{2}}
            \sqrt{ \frac{\log (6c_1/\eta)}{c_2 N}}
            \right\}
            ,
        \end{align}
        so in total with probability at least $1-\frac{2\eta}{3}$
        \begin{align}\notag
            S_2(\lambda, Z_{1:N})
            & \leq
            \frac{1}{\lambda}
            \max\left\{
            \left(\frac{6c_1}{\eta N^{m-1}}
            (2{\mc{B}(\lambda)}^{1/2})^m C_m
            \right)^{1/m},
            \sqrt{{\mc{B}(\lambda)} C_{2}}
            \sqrt{ \frac{\log (6c_1/\eta)}{c_2 N}}
            \right\}^2
            \\  
            & \lesssim 
            \frac{{\mc{B}(\lambda)}}{\lambda}
            \max\left\{
            \frac{1}{(\eta N^{m-1})^{1/m}}
            ,
            \sqrt{ \frac{\log (6c_1/\eta)}{N}}
            \right\}^2
            .
        \end{align}
        We now deal with $S_1$. As outlined in the proof of Theorem \ref{theorem:excess-risk:Bernstein}, the first term on the right hand side of \eqref{eq:S1ineq} is bounded by $4$ under the conditions for $S_2$, which hold with probability at least $1-\frac{2}{3}\eta$.
        
        The second term in \eqref{eq:S1ineq} equals $\left\|\frac{1}{N}\sum_{i=1}^N \zeta_3^i \right\|_{\mch}$, where $\zeta_3$ is defined as in the proof of Theorem \ref{theorem:excess-risk:Bernstein} and we recall that $\zeta_3$ is centered.
        Further, due to \eqref{eq:zeta3bound}, we have
        \begin{align}\notag
            \mathbb{E}\big[\|\zeta_3\|_{\mch}^m\big]
            \leq
            \frac{1}{\lambda^{m/2}}
            \int_{\mcx}
            k(\bx,\bx)^{m/2}
            \left(
            \int_{\R}
            |y - f^{\ast}(\bx)|^{m}
            \, \nu(\mathrm{d}y|\bx)
            \right)
            \, \nu_{\mc{X}}(\mathrm{d}\bx)
            \leq \frac{\sqrt{C_{m}} D_m}{\lambda^{m/2}}.
        \end{align}
        Thus, Theorem \ref{theorem:fuknagaev} implies that with probability at least $1-\frac{1}{3}\eta$, we have
        \begin{align}\notag
            \left\|
            (\mc{T}+\lambda)^{-1/2}
            \big(g_{Z_{1:N}} - \mc{T}_{\bX_{1:N}} f^{\ast}\big)
            \right\|_{\mch}^{2} &
            \leq  
            \max\left\{
            \left(\frac{6c_1 }{\eta N^{m-1}} \frac{\sqrt{C_{m}} D_m}{\lambda^{m/2}}\right)^{1/m},
            \sqrt{\frac{C_{1} D_2}{\lambda}}
            \sqrt{ \frac{\log (6c_1/\eta)}{c_2 N}}
            \right\}^2
            \\ \notag & 
            \lesssim 
            \frac{1}{\lambda}
            \max\left\{
            \left(\frac{1 }{\eta N^{m-1}}\right)^{1/m},
            \sqrt{ \frac{\log (6c_1/\eta)}{N}}
            \right\}^2
        \end{align}
        Hence, the Theorem follows by combining the previous estimates,  since now
        \begin{align}\notag
            \Erel[f_{N}^\lambda]
            \lesssim & \bigg( \mc{A} (\lambda) +\left( \frac{1}{\lambda}
            +\frac{{\mc{B}(\lambda)}}{\lambda}\right)
            \max\left\{
            \left(\frac{1 }{\eta N^{m-1}}\right)^{1/m},
            \sqrt{ \frac{\log (6c_1/\eta)}{N}}
            \right\}^2 
            \bigg).
        \end{align}
    \end{proof}
    \begin{proof}[Proof of Theorem \ref{theorem:rates:heavytail}]
        As in the proof of Theorem \ref{theorem:rates:Bernstein}, 
        plugging the bound \eqref{eq_bounds_A_B} into Theorem \ref{theorem:errorbound:heavytail} yields that with probability at least $1-\eta$ we have
        \begin{align}\notag
            \Erel[f_{N}^\lambda]
            \lesssim & \left( \lambda \|f^{\ast}\|_{\mch}^2 +\left( \frac{1}{\lambda}
            +\frac{\|f^{\ast}\|_{\mch}^2}{\lambda}\right)
            \max\left\{
            \left(\frac{1 }{\eta N^{m-1}}\right)^{1/m},
            \sqrt{ \frac{\log (6c_1/\eta)}{N}}
            \right\}^2 
            \right),
        \end{align}
        provided $\lambda,N$ fulfill the conditions of Theorem \ref{theorem:errorbound:heavytail}. Plugging the choices of $\lambda_N, N$ proposed in Theorem \ref{theorem:rates:heavytail} into the above bound and noting that they fulfill the condition of Theorem \ref{theorem:errorbound:heavytail}, the result follows, since
        \begin{align}
            \Erel[f_{N}^{\lambda_N}]
            \lesssim & \lambda_N \|f^{\ast}\|_{\mch}^2 +
            \frac{1}{\lambda_N}
            \left( 1
            +\|f^{\ast}\|_{\mch}^2\right)
            \frac{\log (6c_1/\eta)}{N}
            \\ \notag 
            \lesssim & 
            \frac{\sqrt{\log (6c_1/\eta)}}{\sqrt{N}}
            \bigg( 
            \|f^{\ast}\|_{\mch}^2 
            \left(
            C''_m
            \vee
            \left( \frac{1}{\|f^{\ast}\|_{\mch}^2}
            +1\right)^{1/2}
            \right)
            \\ \notag
            & +
            \frac{1}{C_m''}
            \left( 1
            +\|f^{\ast}\|_{\mch}^2\right)
            \bigg).
        \end{align}
    \end{proof}
    We close this section with the proof of Theorem \ref{theorem:consistency:heavytail}.
    \begin{proof}[Proof of Theorem \ref{theorem:consistency:heavytail}]
        To obtain an analogous bound to that of Theorem \ref{theorem:errorbound:heavytail}, we need to employ a moment bound on the random variable $\zeta_2 = \mc{T}_{\bX} (f^\lambda- f^{\ast})$ in the proof of Theorem \ref{theorem:errorbound:heavytail}, respectively Theorem \ref{theorem:excess-risk:Bernstein}, that is in terms of $L^p(\nu_\mcx)$ norms instead of $\mch$ norms. To this end, observe that we have
        \begin{align}
            & \mb{E}\big[\|\zeta_2-\mb{E}[\zeta_2]\|^3_{\mch}\big]
            =
            \mb{E}\big[\|(\mc{T}_{\bX}-\mc{T}) (f^\lambda-f^{\ast})\|^3_{\mch}\big]
            \\ \notag \leq &
            2^3 \mb{E}\big[k (\bX,\bX)^{3/2} | f^\lambda (\bX)-f^{\ast} (\bX)|^3\big]
            \leq 
            \sqrt{4^3 C_3} \,\|f^\lambda -f^{\ast}\|_{L^6(\nu_\mcx)}^3
            , 
        \end{align}
        so by the Fuk-Nagaev inequality, i.e., Theorem \ref{theorem:fuknagaev}, for $0< \lambda \leq \lambda_0$ we have
        \begin{align}\label{eq:consistency:adjustment}
            \| (\mc{T}-\mc{T}_{\bX_{1:N}})(f^\lambda-f^{\ast}) \|_{\mch}
            \lesssim 
            \max
            \left\{
            \left(\frac{1}{\eta^{1/3} N^{2/3}}\right),
            \sqrt{ \frac{\log (6c_1/\eta)}{N}}
            \right\},
        \end{align}
        which scales as $1/\sqrt{N}$ for large enough $N$, as the previous bound \eqref{eq:S2:boundpart2:heavy}.
        In particular, we still obtain an adjusted version of the bound of Theorem \ref{theorem:excess-risk:Bernstein}, just with different constants in the last term of the right hand side of Theorem \ref{theorem:errorbound:heavytail}'s bound, which is irrelevant for asymptotic consistency, and without the multiplicative $\mc{B}(\lambda)$ term.
        
        Since by Lemma \ref{lemma_A_B_bound} we have $\lim_{\lambda\to 0}\mc{A}(\lambda)=0$, the right hand side of the adjusted version of Theorem \ref{theorem:errorbound:heavytail}'s bound converges to zero for any sequence $\{\lambda_N\}_{N\in \mb{N}}\subset (0,\infty)$ fulfilling $\lim_{N\to \infty}\lambda_N N=\infty$ and for which \eqref{condition:errorbound:heavytail} holds for any $\eta\in (0,1)$ for large enough $N$. 
        
        Note that a sufficient condition for \eqref{condition:errorbound:heavytail} to hold for $N$ large enough is $\lim_{N\to \infty}\lambda_N^2 N=\infty$ and $\lim_{N\to \infty}1/\lambda_N=\infty$, so we can replace the former asymptotic conditions \eqref{condition:errorbound:heavytail} by the latter, i.e., \eqref{eq_condition_no_rates:heavytail} suffices for consistency.
    \end{proof}

\subsubsection{Signature kernel moments for Brownian data}\label{sec:sig-kernel-moment:BM}
We prove Lemma \ref{lemma:BM-kernel-moments}.
\begin{proof}[Proof of Lemma \ref{lemma:BM-kernel-moments}]
We make use of the concepts of words and shuffle products of words, denoted for words $v,w\in W$ by $v \shuffle w$. 
We refer the reader unfamiliar with the concepts of words and shuffles to
\cite[Section 1.1]{BayerPelizzariSchoenmakers2025}.
Throughout this proof, we write $W$ for the set of all words, \textit{excluding} the empty word $w = \emptyset$. For $n\in \N$, we write $W_n$ for the subset containing all words of length $n$. Given a collection of words $w_1,\dots,w_m$, we will denote the multi-shuffle of the words by $\textnormal{Sh}(w_1,\dots,w_m) \coloneqq w_1 \shuffle\dots \shuffle w_m$, which is well defined by associativity.

Denoting by $\bar{B}$ a time augmented Brownian motion taking values in $\mb{R}^{d+1}$ with time horizon $[0,T]$, for $m\in \mb{N}$ we aim to bound $\mb{E}\left[k(\bar{B},\bar{B})^m\right]$.
Since all estimates in this proof hold when replacing $T$ by $T\vee 1$ in the upper bounds, we will assume for clarity of exposition $T\geq 1$. Further, for the same reason, setting $k^- (\cdot,\cdot) \coloneqq k(\cdot,\cdot)-1$, we will use that
\begin{equation}\label{eq:bound-k-one-k}
    \mb{E}[k(\bar{B},\bar{B})^m] \leq 2^{m-1}(1+ \mb{E}[k^-(\bar{B},\bar{B})^m])
\end{equation}
and focus on bounding $\mb{E}[k^-(\bar{B},\bar{B})^m]$.
To this end, note first that we can expand
\begin{align}\label{eq_expansion_strat_sigkernel}
    \mb{E}[k^-(\bar{B},\bar{B})^m]
    = \sum_{w_1\in W} \dots \sum_{w_m\in W} 
    \mb{E}\left[
    \prod_{i\leq m} \sig^{w_i}(\bar{B})^2
    \right].
\end{align}
Iteratively applying the shuffle product, we find
\begin{align}\label{eq_multi_shuffle_exp_strat}
    \mb{E}\left[
    \prod_{i\leq m} \sig^{w_i}(\bar{B})^2
    \right]
    = \sum_{w^\ast \in \textnormal{Sh}(w_1,w_1,\dots,w_m,w_m)}
    \mb{E}\left[
    \sig^{w^\ast}(\bar{B})
    \right],
\end{align}
where $\textnormal{Sh}(w_1,w_1,\dots,w_m,w_m)$ denotes the set of multi-shuffles of the $2m$ words $w_1,w_1,\dots,w_m,w_m$.
In the following, we will write $\textnormal{MS}_m\coloneqq\textnormal{Sh}(w_1,w_1,\dots,w_m,w_m)$, making the dependence of the set of multi-shuffles on the words implicit.

We want to compute the left hand side of \eqref{eq_expansion_strat_sigkernel} using the above identity \eqref{eq_multi_shuffle_exp_strat}. This reduces the problem to that of computing expectations of Stratonovich signatures for time augmented Brownian motion. The standard way to do so is applying the Fawcett formula \cite{Fawcett2003, LyonsVictoir2004}. 
Recall that Fawcett states that the expectation of the Stratonovich signature of time augmented Brownian motion is given by a tensor exponential, namely that
\begin{align}\label{eq_fawcett}
\mathbb{E}\!\left[\operatorname{Sig}(\bar{B})\right]
=
\exp\!\left(
T e_{d+1}
+
\frac{T}{2}\sum_{i=1}^{d} e_i \otimes e_i
\right).
\end{align}
Hence, the entries of the expected signature are non-vanishing only for words that consist of pairs of letters corresponding to Brownian motion and letters corresponding to time. The associated coefficient for a word of length $L$ with $n$ pairs of space indices and $L-2n$ occurrences of the time index is by Fawcett
\begin{align}
    \tilde{c}_{n,L} \coloneqq \frac{1}{(L-n)!} \frac{T^{L-n}}{2^n}
    = 
    \frac{1}{(L-n)!} \frac{T^{L}}{(2T)^n},
\end{align}
where $0\leq n \leq L/2$ is the number of pairs.
Assuming w.l.o.g. $T\geq 1$, an upper bound on this quantity is given by
\begin{align}
    \tilde{c}_{n,L} \leq \frac{T^{L}}{(L/2)!} \eqqcolon c_L.
\end{align}
Assume that we consider words $w_1,\dots,w_m$ with total length $L=l_1+\dots + l_m$ such that all words consist only of the letter $d+1$ corresponding to time. 
For such words, all elements of the multi-shuffle are equal to the word with all letters equal $d+1$ with doubled length and thus associated coefficient positive and equal $\frac{T^{2L}}{(2L)!}$.
Hence, one finds that 
\begin{align}
    \mb{E}\left[
    \prod_{i\leq m} \sig^{w_i}(\bar{B})^2
    \right]
    = \# \textnormal{MS}_m \frac{T^{2L}}{(2L)!}.
\end{align}
Since by computation also 
\begin{align}
    \mb{E}\left[
    \prod_{i\leq m} \sig^{w_i}(\bar{B})^2
    \right]
    =
    \prod_{i\leq m} 
    \mb{E}\left[
    \sig^{w_i}(\bar{B})
    \right]^2
    = \prod_{i=1}^m
    \frac{T^{2l_i}}{(l_i!)^2},
\end{align}
we have in particular\footnote{Observe that one can derive $\# \textnormal{MS}_m$ also by a combinatorial argument not involving any probability.}
\begin{align}
    &
    \#\textnormal{MS}_m 
    =
    (2L)!
    \prod_{i=1}^m
    \frac{1}{(l_i!)^2}
    \leq 
    (2L)!
    \prod_{i=1}^m
    \frac{1}{((L/m)!)^2}
    \\ \notag =&
    (2L)!
    \frac{1}{((L/m)!)^{2m}}.
\end{align}
In the above we used that for $l_i$ fulfilling $l_1+\dots+l_m=L$, we have by log-convexity of $c \mapsto c!$ that
\begin{align}
    \prod_{i=1}^m (L/m)!
    \leq \prod_{i=1}^m l_i!.
\end{align}
Since $c_{2L}$ bounds above $\mb{E}\left[\sig^{w}(\bar{B})\right]$ if $w$ is a word of length $2L$ and since   $\#\textnormal{MS}_m$ does not depend on the choice of words but just on their lengths, we obtain from this that for any words $w_1,w_2,\dots, w_m$ that are of total length $L=l_1+\dots+l_m$, we have
\begin{align}
    \mb{E}\left[
    \prod_{i\leq m} \sig^{w_i}(\bar{B})^2
    \right]
    \leq &
    (2L)!
    \frac{1}{((L/m)!)^{2m}}
    c_{2L}
    \\ \notag = &
    (2L)!
    \frac{1}{((L/m)!)^{2m}}
    \frac{T^{2L}}{L!}
    \\ \notag = &
    \frac{(2L)!}{L! ((L/m)!)^{2m} }\,T^{2L}.
\end{align}
Combining the above results and writing $L=l_1+\dots+l_m$, we have
\begin{align}\notag
    \mb{E}[k^-(\bar{B},\bar{B})^m]
    =&
    \sum_{w_1\in W} \dots \sum_{w_m\in W} 
    \mb{E}\left[
    \prod_{i\leq m} \sig^{w_i}(\bar{B})^2
    \right]
    \\ \notag \leq &
    \sum_{w_1\in W} \dots \sum_{w_m\in W}
    \frac{(2L)!}{L! ((L/m)!)^{2m}}\,T^{2L}.
\end{align}
To bound the above right hand side, note that the number of collections of words $w_1,w_2,\dots, w_m$ such that $L=l_1+\dots+l_m$ is given by $(d+1)^L \binom{L-1}{m-1}$, since each such collection of words can be obtained by splitting a word of length $L$ into $m$ subwords by choosing $m-1$ splitting letters.
In particular, one thus finds that
\begin{align}\notag
    \mb{E}[k^-(\bar{B},\bar{B})^m]
    \leq &
    \sum_{L\in \mb{N} \setminus {0}}
    (d+1)^L \binom{L-1}{m-1}
    \frac{(2L)!}{L! ((L/m)!)^{2m}}\,T^{2L}
    \\ \notag \leq &
    \sum_{L\in \mb{N} \setminus {0}}
    ((d+1)T^2)^L \frac{L^m}{(m-1)!}
    \frac{(2L)!}{L! ((L/m)!)^{2m}}
    \\ \label{eq_C_bound_stirling} \lesssim &
    \sum_{L\in \mb{N}\setminus \{0\}}
    ((d+1)T^2)^L \frac{L^m}{(m-1)!}
    \frac{\sqrt{2}m^m}{(2\pi L)^m}
    \cdot \left(
    \frac{4m^2 e}{L}
    \right)^L.
\end{align}
Indeed, writing $L=n m$ with $n\in \frac{1}{m}\cdot \N$, by Stirling's formula we have
\begin{align}
    & \frac{(2L)!}{L! ((L/m)!)^{2m}}
    =
    \frac{(2n m)!}{(n m)!   (n!)^{2m}}
    \\ \notag \lesssim & 
    \sqrt{2\pi 2n m} \left( \frac{2n m}{e} \right)^{2n m}
    \frac{1}{\sqrt{2\pi n m}}
    \left( \frac{e}{n m} \right)^{n m}
    \frac{1}{(2\pi n)^m}
    \left( \frac{e}{n} \right)^{2nm}
    \\ \notag = & 
    \sqrt{2} 
    \left( 2 m \right)^{2n m}
    \left( \frac{e}{n m} \right)^{n m}
    \frac{1}{(2\pi n)^m}
    \\ \notag = & 
    \sqrt{2} 
    \left( 4m^2 \right)^{n m}
    \left( \frac{e}{n m} \right)^{n m}
    \frac{1}{(2\pi n)^m}.
\end{align}
Substituting $n=L/ m$, we obtain 
\begin{align}
    \frac{(2L)!}{L! ((L/m)!)^{2m}}
    \lesssim 
    &
    \sqrt{2} 
    \left( 4m^2 \right)^{n m}
    \left( \frac{e}{n m} \right)^{n m}
    \frac{1}{(2\pi n)^m}
    \\ \notag
    =
    &
    \sqrt{2} 
    \left( 4m^2 \right)^{L}
    \left( \frac{e}{L} \right)^{L}
    \frac{1}{(2\pi L/m)^m}
    \\ \notag
    =
    &
    \frac{\sqrt{2}m^m}{(2\pi L)^m}
    \cdot \left(
    \frac{4 m^2 e}{L}
    \right)^L,
\end{align}
where all constants hidden by $\lesssim$ are universal, coming from the Stirling formula, so indeed \eqref{eq_C_bound_stirling} holds. 

Applying Stirling once more on $m! = m(m-1)!$, e.g. using $m^{m}/m! \sim e^m/\sqrt{2\pi m}$, we find that
\begin{align}
    \mb{E}[k^-(\bar{B},\bar{B})^m]
    \lesssim & 
    \sum_{L\in \mb{N}\setminus \{0\}}
    e^m \frac{\sqrt{m}}{\sqrt{2\pi}}
    ((d+1)T^2)^L 
    \frac{\sqrt{2}L^m}{(2\pi L)^m}
    \cdot \left(
    \frac{4m^2 e}{L}
    \right)^L
    \\ =& 
    e^m \frac{\sqrt{m}}{\sqrt{2\pi}}
    \frac{\sqrt{2}}{(2\pi)^m}
    \sum_{L\in \mb{N}\setminus \{0\}}
    \left(
    \frac{4m^2 e (d+1)T^2}{L}
    \right)^L.
\end{align}
Using that for $c\geq 0$ we have 
\begin{align}
    \sum_{L\in \mb{N}\setminus \{0\}}
    \left(
    \frac{c}{L}
    \right)^L
    \leq 
    \sum_{L\in \mb{N}\setminus \{0\}}
    \frac{c^L}{L!}
    =
    e^c -1,
\end{align}
the above implies 
\begin{align}
    & \mb{E}[k^-(\bar{B},\bar{B})^m]
    \lesssim 
    \frac{\sqrt{m}}{(2\pi)^m}
    \exp(4m^2 e (d+1)T^2+ m).
\end{align}
Combined with \eqref{eq:bound-k-one-k}, we obtain
\begin{equation*}
    \mb{E}[k(\bar{B},\bar{B})^m]
    \leq 2^{m-1}(1+ \mb{E}[k^-(\bar{B},\bar{B})^m])
    \lesssim 
    2^{m-1}
    \left(1+
    \frac{\sqrt{m}}{(2\pi)^m}
    \exp(4m^2 e (d+1)T^2+ m)\right),
\end{equation*}
which is the claim of the Lemma.
    
\end{proof}

\subsubsection{Robust signature kernels and proofs}\label{appendix:robust}
In this section, we give an overview of robust
signatures and prove Theorem \ref{theorem:rates:robust} for kernel ridge regression with robust signatures. 

The idea of robust signatures is to construct a normalized feature map for path-valued data based on signatures, preserving their algebraic and geometric properties.
While the signature is not statistically robust w.r.t. outliers, the robust signature is so in a precise statistical sense \cite{chevyrev_oberhauser_robust}, and maintains at the same time the approximation quality of signatures \cite{chevyrev_oberhauser_robust, BayerPelizzariSchoenmakers2025}.

To construct the robust signature, one does not normalize the signature in feature space directly, but considers instead the signature of a suitably normalized path.
Given a fixed $\kappa>1$,
the normalization is obtained by identifying a function $\blambda:T((\mathbb{R}^{d+1}))\to \mb{R}_+$ that is such that the mapping
\begin{align}
    \widehat{\sig}: \mcx \to B_{\sqrt{\kappa}}\left(T((\mb{R}^{d+1}))\right),
    \quad
    \widehat{\sig}:
    \bx\mapsto 
    \sig\left(\blambda\left(\sig(\bx)\right)\cdot \bx\right)
\end{align}
is injective up to equal signatures, i.e., $\widehat{\sig}(\bx)=\widehat{\sig}(\bx')$ if and only if $\sig(\bx)=\sig(\bx')$. 
Then, $\widehat{\sig}\left(\bx\right)$ is the signature of the rescaled path $\blambda\left(\sig(\bx)\right)\cdot \bx$ \footnote{Observe that $\blambda$ considered in this context is unrelated to the ridge regression parameter $\lambda$. Since we will not deal explicitly with $\blambda$ but only with the induced robust signatures, we will stick to this classical notation.}.
We make here and in the following implicit the dependence on $\kappa>1$.
Note that one can also identify the robust signature as a tensor normalization of the signature, since
\begin{align}
    \left\|\sig\left(\blambda(\sig(\bx))\,\bx\right)\right\|^2_{T^2((\mb{R}^{d+1}))}
    = \sum_{k\in \mb{N}}
    \blambda(\sig(\bx))^{2k}
    \|\sig(\bx)^{(k)}\|^2_{(\mb{R}^{d+1})^{\otimes k}}.
\end{align}
As in the classical case, the robust signature kernel is obtained as the scalar product of robust signatures, i.e., for paths $\bx,\bx'$ we set
\begin{align}
    \hat{k} (\bx,\bx') \coloneqq \left\langle \widehat{\sig}(\bx),
    \widehat{\sig}(\bx') \right\rangle_{T^2((\mb{R}^{d+1}))}.
\end{align}
Note that different to the signature kernel, the above kernel is bounded, i.e., $\hat{k} (\bx,\bx) \leq \kappa$.

In practice, one considers the robust signature kernel truncated at level $m\in \mathbb{N}$, i.e., only the first $m$ levels of the signature are used, so
\begin{align}
    \hat{k}_{m} (\bx,\bx') = \sum_{k\leq m} 
    \left\langle\widehat{\sig}(\bx)^{(k)},
    \widehat{\sig}(\bx')^{(k)} 
    \right\rangle_{(\mb{R}^{d+1})^{\otimes k }}.
\end{align}
The truncated robust signature and hence the truncated robust signature kernel can be computed efficiently using the approach proposed in \cite{chevyrev_oberhauser_robust}, see also \cite{bayer2026local}. The approach relies on identifying $\blambda$ as the root of a polynomial. Note that the Goursat-PDE approach for the signature kernel \cite{goursat_PDE} is not applicable directly to the robust kernel, as the normalization factor $\blambda(\sig(\bx))$ is unknown a priori.

The setup of signature kernel ridge regression translates directly to the robust case, i.e., we aim at minimizing the optimization objective \eqref{eq:empiricallossridge}, but across the RKHS of the robust signature kernel. In general, all objects are defined analogously by replacing the signature kernel with the robust signature kernel everywhere. As mentioned before, we indicate the analogues for the robust kernel by a hat accent, i.e., the analogues of $f^\lambda,f_{N}^\lambda$ are denoted by $\hat{f}^\lambda,\hat{f}_{N}^\lambda$.

Most of the results discussed for classical signatures translate directly up to minor changes to the robust case. 
We start with the following, dealing with the structure of the RKHS of the robust signature kernel.
\begin{lemma}\label{lemma_RKHS_robust}
    For $\alpha \in (\frac{1}{3},1]$, $T>0$,
    consider the stopped rough-path space $\mcx=\Lambda_T^\alpha$ or the time augmented rough path space $\mcx =
    \widehat{\mathscr{C}}^{\alpha}_g([0,T],\R^{d+1})$. 
    Then, the RKHS $\widehat{\mch}$ associated to the robust signature kernel $\hat{k}$ is given as the closure of the function space
    \begin{equation*}
        \widetilde{\mch} \coloneqq \textnormal{span} \{\hat{k}_{\bx} :\bx'\mapsto \hat{k}(\bx,\bx'): \quad \bx \in \mcx\},
    \end{equation*}
    equipped with the scalar product $\langle \cdot, \cdot \rangle_{\widehat{\mch}}$, defined for $\hat{k}_{\bx},\hat{k}_{\bx'}$ by
    \begin{equation*}
        \langle \hat{k}_\bx, \hat{k}_{\bx'}\rangle_{\widehat{\mch}} \coloneqq \hat{k}(\bx,\bx')
    \end{equation*}
    and extended to $\widetilde{\mch}$ by linearity.
    
    Further, we have 
    \begin{equation*}
        \widehat{\mch}= \big\{ f:\mcx \to \R: \quad \,f(\bx) = \langle \ell_f,\sigrob(\bx)\rangle_{T^{2}((\R^{d+1}))} 
        ,\, \ell_f\in T^{2}((\R^{d+1}))
        \big\},
    \end{equation*}
    
    and, given $f\in \widehat{\mch}$, that
    \begin{equation*}
        \|f\|_{\widehat{\mch}} = \inf 
        \left\{\|\ell_f\|_{T^{2}((\R^{d+1}))}:\quad 
        f(\cdot)=\langle \ell_f,\sigrob (\cdot) \rangle_{T^{2}((\R^{d+1}))}, \,
        \ell_f \in T^{2}((\R^{d+1}))
        \right\}.
    \end{equation*}
\end{lemma}
As for Lemma \ref{lemma_RKHS}, the result follows by \cite[Theorem 4.21]{steinwart2008support}. Note that, different to $\mch$ from Lemma \ref{lemma_RKHS}, the elements of $\mchrob$ are bounded functions.

We next state the $L^p$-approximation result for robust signatures due to Bayer et al. \cite{BayerPelizzariSchoenmakers2025}. It requires the following assumption, which holds in all practically relevant cases, in particular for Brownian motion.
\begin{assumption}[{\cite[Assumption 2.6]{BayerPelizzariSchoenmakers2025}}]\label{assumption_bayer} The following two conditions hold.
    \begin{enumerate}
    \item $\alpha \in (1/3,1)$ is such that
    \[
    q := \frac{1}{\alpha} \notin \mathbb{N}.
    \]
    
    \item $\mu$ is a measure on the Borel space of stopped rough paths
    $\left(\Lambda_T^\alpha, \mathcal{B}(\Lambda_T^\alpha)\right)$ such that for all $\alpha < \beta < \frac{1}{\lfloor q \rfloor}$ we have
    \[
    \mu(\Lambda_T^\alpha) < \infty
    \quad \text{and} \quad
    \mu\bigl(\Lambda_T^\alpha \setminus \Lambda_T^\beta \bigr) = 0.
    \]
    \end{enumerate}
\end{assumption}
The approximation result for robust signatures then reads as follows.
\begin{theorem}[{\cite[Theorem 2.7]{BayerPelizzariSchoenmakers2025}}]
    Let $\alpha \in (1/3,1)$ and let $\mu$ be a probability measure on
    $\bigl(\Lambda_T^\alpha, \mathcal{B}(\Lambda_T^\alpha)\bigr)$ such that
    Assumption \ref{assumption_bayer} holds true. Then, for all
    $f \in L^{p}(\Lambda_T^{\alpha}, \mu)$ with $1 \le p < \infty$, there exists
    a sequence $(f_n)_{n \in \mathbb{N}} =(\langle \ell_{f_n}, \sigrob (\cdot)\rangle)_{n \in \mathbb{N}} \subset \mchrob$ with $\ell_{f_n}\in T(\R^{d+1})$,
    such that
    \[
    \lVert f - f_n \rVert_{L^{p}(\Lambda_T^{\alpha}, \mu)} \to 0 \quad \text{as } n \to \infty.
    \]
\end{theorem}
The proof relies on a weighted version of the Stone-Weierstrass result, see Theorem 26 in \cite{chevyrev_oberhauser_robust}.

We close this section with the proof of Theorem 
\ref{theorem:rates:robust}. The central ingredient is the following.
\begin{theorem}[{\cite[Theorem 4]{Caponnetto2007OptimalRates}}]\label{theorem:excess-risk:robust}
    Let $N \in \mathbb{N}$, $\lambda > 0$, $\eta  \in (0,1]$,
    $f^{\ast}\in \widehat{\mch}$, and assume that
     Assumption \ref{assumption_nu_finite_sample:robust} holds. Then, with probability at least $1 - \eta$, we have
    \begin{align}\label{eq_risk_bound_robust}
    \mathcal{E}[\hat{f}_{N}^{\lambda}] - \mathcal{E}[f^{\ast}] \leq C_{\eta} 
    \left(
    \mc{A} (\lambda) + \frac{\kappa^{2} \mc{B}(\lambda)}{N^{2} \lambda} + \frac{\kappa \mc{A}(\lambda)}{N \lambda} + \frac{\kappa M^{2}}{N^{2} \lambda} + \frac{\Sigma^2 \mathcal{N}(\lambda)}{N}
    \right)
    \end{align}
    provided that
    \begin{equation}\label{eq:condition:error:robust}
        N \geq \frac{2 C_{\eta} \kappa \mathcal{N}(\lambda)}{\lambda} \quad \text{and} \quad \lambda \leq \|\widehat{\mc{T}}\|_{\mathcal{L}(\mchrob)},
    \end{equation}
    where $C_{\eta} = 96 \log^{2}(6/\eta)$.
    Here, $\mc{A}(\lambda),\mc{B}(\lambda)$ are defined as in \eqref{eq:residual}, \eqref{eq_reconstruction_error}, whereas the effective dimension $\mc{N}(\lambda)$ is defined as 
    \begin{align}
        \mc{N}(\lambda) \coloneqq \textnormal{Tr}\left(
        (\widehat{\mc{T}}+\lambda)^{-1} \widehat{\mc{T}}
        \right).
    \end{align}
\end{theorem}
\begin{remark}
    We use the loose bound $\mc{N}(\lambda) \leq \frac{\textnormal{Tr}(\widehat{\mc{T}})}{\lambda}$ for the effective dimension. In \cite{Caponnetto2007OptimalRates}, sharper bounds are derived under assumptions on the spectrum of the operator $\widehat{\mc{T}}$.
\end{remark}
  We now give the proof of Theorem \ref{theorem:rates:robust}.
\begin{proof}[Proof of Theorem \ref{theorem:rates:robust}]
    Since the choice \eqref{eq:lambdaN:robust} of $\lambda=\lambda_N$ and condition \eqref{eq:conditionN:robust} ensure that \eqref{eq:condition:error:robust} holds, Theorem \ref{theorem:excess-risk:robust} is applicable.
    We thus start computing the bound \eqref{eq_risk_bound_robust}. Recall that in Lemma \ref{lemma_A_B_bound} we established already bounds on $\mc{A}(\lambda),\mc{B}(\lambda)$. For $\mc{N}(\lambda)$, we have
    \begin{align}\label{eq_bound_effective_dim}
        \mc{N}(\lambda) = 
        \textnormal{Tr}\left(
        (\widehat{\mc{T}}+\lambda)^{-1} \widehat{\mc{T}}
        \right)
        \leq \textnormal{Tr}\left(
        (\lambda)^{-1} \widehat{\mc{T}}
        \right)
         = \frac{\textnormal{Tr}(\widehat{\mc{T}})}{\lambda}.
    \end{align}
    Plugging this bound together with the bounds of Lemma \ref{lemma_A_B_bound} in \eqref{eq_risk_bound_robust}, we find that
    \begin{align}\notag
        \Erel [\hat{f}_N^{\lambda}] \leq C_{\eta} 
        \left(
        \lambda \|f^{\ast}\|_{\mchrob}^2 + \frac{\kappa^{2} \|f^{\ast}\|_{\mchrob}^2}{N^{2} \lambda} + \frac{\kappa \|f^{\ast}\|_{\mchrob}^2}{N} + \frac{\kappa M^{2}}{N^{2} \lambda} + \frac{\Sigma^2 \textnormal{Tr}(\widehat{\mc{T}})}{N \lambda}
        \right)
    \end{align}
    Since the second condition in \eqref{eq:condition:error:robust} implies $\mc{N}(\lambda)\geq \frac{1}{2}$ and since $C_\eta \geq 96 \log^2(6)\geq 102$, the first condition implies
    \begin{align}\label{eq:condition:lambdaN:clean}
        \lambda N \geq C_\eta \kappa\geq 102 \kappa.
    \end{align}
    Combining this with the inequalities $\lambda \leq \|\widehat{\mc{T}}\|_{\mathcal{L}(\mchrob)}\leq \Tr (\widehat{\mc{T}}) \leq \kappa$, we find that
    \begin{align}
        \Erel [\hat{f}_N^{\lambda}]
        \leq 
        \frac{6}{5}C_\eta\left( \lambda \|f^{\ast}\|_{\mchrob}^2 
        +
        (\frac{M^2}{102}+\Sigma^2)
        \frac{\Tr (\widehat{\mc{T}})}{N \lambda}
        \right).
    \end{align}
    Plugging the formula for $\lambda=\lambda_N$ from Theorem \ref{theorem:rates:robust} into the above and simplifying, the result follows. 
\end{proof}

\subsection{Mercer representations for signature kernels}\label{section_sub_mercer}
We discuss how the work \cite{SteinwartScovel2012} implies beyond Lemma \ref{lemma:steinwart-scovel} Mercer-type series expansions for the signature kernel.
Consider a compact metric space $(X,d)$, a continuous kernel $k:X\times X \to \R$ with associated separable RKHS $\mch$, and a finite strictly positive Borel measure $\nu$. Then, the integral operator $\mc{T}:L^2(\nu_\mcx)\to \mch$ defined by 
\begin{align}\label{eq_def_T_k}
    \mc{T} f\coloneqq \int_X k(\cdot,x') f(x') \,\nu(\mathrm{d}x') \in C(X) \cap \mch
\end{align}
is positive, compact, and self-adjoint. 
By the spectral Theorem, i.e., Lemma \ref{lemma:steinwart-scovel}, we obtain an expansion for $\mc{T}$ in terms of eigenvectors and eigenvalues $([e_i]_{\sim})_{i\in \mb{N}},\,(\mu_i)_{i\in \mb{N}}$ that is of the form \eqref{eq_T_k_L2_classical_2}, where $(e_i)_{i\in \mb{N}}\subset \mch\subset C(X)$.
The Mercer Theorem, see also \cite{RieszNagy1990}, states that one can express the kernel $k$ as a series, using the continuous representatives and associated eigenvalues of the spectral decomposition for $\mc{T}$. More precisely, for all $x,x'\in X$ one has
\begin{align}\label{eq_Mercer_classical}
    k(x,x') = \sum_{i\in \mb{N}} \mu_i e_i(x) e_i(x'),
\end{align}
where the sum converges absolutely and uniformly in $x,x'\in X$. 
The Mercer representation \eqref{eq_Mercer_classical} has a number of consequences \cite{SteinwartScovel2012}, for example, one can use it to characterize explicitly the RKHS $\mc{H}$, see e.g. \cite{cucker2002mathematical, CuckerZhou2007}. 

Steinwart and Scovel \cite{SteinwartScovel2012} studied Mercer representations in full generality without assuming the space $X$ to be compact. Hence, their results imply Mercer representations for signature kernels on (stopped) rough path space.
This is the content of the next Lemma. 
\begin{lemma}\label{lemma_Mercer_sig}
    Let $\mc{X}=\Lambda^{\alpha}_T$ or $\mc{X}=\widehat{\mathscr{C}}^{\alpha}_g([0,T],\mb{R}^{d+1})$ with $\alpha \in (\frac{1}{3},1]$, $T>0$, $\nu_{\mc{X}}\in \mc{P}(\mc{X})$, and $k=k_{\sig,\phi}$ be the weighted signature kernel on $\mc{X}$, where $\phi:\N\to \R_{+}\setminus \{0\}$ fulfills \eqref{eq_condition_phi}. Assume that \eqref{ass_kernel_diagonal_eq} holds, so that the RKHS $\mc{H}=\mch_\phi$ is compactly embedded into $L^2(\nu_{\mc{X}})$. Furthermore, let $(e_i)_{i\in \mb{N}}\subset \mc{H}$ and $(\mu_i)_{i\in \mb{N}}\subset (0,\infty)$ be as in Lemma \ref{lemma:steinwart-scovel}.
    
    The weighted signature kernel admits a globally convergent Mercer expansion, i.e., for all $\mathbf{x},\mathbf{x}'\in \mc{X}$ we have
    \begin{align}\label{eq_Mercer_theorem}
        k(\mathbf{x},\mathbf{x}')=\sum_{i\in \mb{N}} \mu_i e_i(\mathbf{x})e_i(\mathbf{x}'),
    \end{align}
    if and only if the embedding $\mc{H} \hookrightarrow L^2(\nu_{\mc{X}})$ is injective. In \eqref{eq_Mercer_theorem}, the convergence is absolute and uniform in $\mathbf{x},\mathbf{x}'\in A$ for compact $A \subset \mc{X}$. 
    
    In the above setup, injectivity holds if and only if there exists no $\ell\in T^{2,\phi}((\mb{R}^{d+1}))$ such that $\mathbf{x}\mapsto \langle \ell, \sig(\mathbf{x})\rangle_{T^{2,\phi}((\mb{R}^{d+1}))}\in \mc{H}\setminus \{0\}$ and for $\bX\sim \nu_{\mc{X}}$ one has almost surely
    \begin{align}\label{eq_condition_Mercer}
        \langle \ell, \sig(\bX)\rangle_{T^{2,\phi}((\mb{R}^{d+1}))} = 0.
    \end{align}
\end{lemma}
\begin{remark}\label{remark_mercer}
    One can check easily that a sufficient condition for the injectivity is that $\textnormal{supp}(\nu_{\mc{X}})$ contains an open set. This one can show by a perturbation argument, using that polynomials on open sets can only be constant if they are trivial.
\end{remark}
\begin{proof}[Proof of Lemma \ref{lemma_Mercer_sig}]
    The first part of Lemma \ref{lemma_Mercer_sig} follows by Theorem 3.1 in \cite{SteinwartScovel2012}, which states that a measurable kernel admits a pointwise convergent Mercer representation if and only if the embedding of the RKHS $\mc{H}$ into $L^2(\nu_{\mc{X}})$ is injective. Note that pointwise convergence is here equivalent to absolute convergence and uniform convergence on compacts, see
    also \cite[Corollary 2.10]{SteinwartScovel2012}. 
    
    Since injectivity of a linear map is equivalent to the map having a trivial kernel, it is injective if and only if there exists no $f\in \mc{H}\setminus \{0\}$ such that $f(\bX)=0$ a.s., where $\bX\sim \nu_{\mc{X}}$. Recalling the RKHS representation result contained in Lemma \ref{lemma_RKHS}, this establishes that the non-existence of a functional fulfilling condition \eqref{eq_condition_Mercer} is equivalent to injectivity. 
\end{proof}
The case where injectivity of the embedding $\mc{H} \hookrightarrow L^2(\nu_\mcx)$ is not assumed is treated in the following Lemma.
\begin{lemma}\label{lemma_Mercer_sig_local}
        Consider the setup of Lemma \ref{lemma_Mercer_sig}. Then, for all $\bx\in \textnormal{supp}(\nu_{\mc{X}})$ and $\bx'\in \mc{X}$ we have that
        \begin{align}
            k(\mathbf{x},\mathbf{x}')=\sum_{i\in \mb{N}} \mu_i e_i(\mathbf{x})e_i(\mathbf{x}'),
        \end{align}
        where the convergence is absolute and for compact $A \subset \textnormal{supp}(\nu_{\mc{X}})$, it is uniform if $\bx,\bx'\in A$.
\end{lemma}
\begin{remark}
    Note that in the general setting of Lemma \ref{lemma_Mercer_sig_local}, the restriction to $\textnormal{supp}(\nu_{\mc{X}})\times \mc{X}$ is expected a priori, as the integral operator $\mc{T}$ depends only on the value of $k$ on that set.
\end{remark}
\begin{proof}[Proof of Lemma \ref{lemma_Mercer_sig_local}]
    Since the weighted signature kernel is continuous, the result is a direct consequence of the more general Corollary 3.5 in \cite{SteinwartScovel2012}.
\end{proof}
\newpage
\printbibliography
\end{document}